\documentclass[11pt,twoside, a4paper]{article}

\usepackage{latexsym}
\usepackage[english]{babel}
\usepackage{t1enc}
\usepackage{mathrsfs}
\usepackage{ifthen}
\usepackage{url}
\usepackage{epsfig}
\usepackage{amsbsy}
\usepackage{extarrows}
\usepackage{amsfonts}
\usepackage{amsmath}
\usepackage{amssymb}
\usepackage{amsthm}
\usepackage{amsxtra}
\usepackage{bbm}
\usepackage{color}
\usepackage{relsize} 
\usepackage{enumitem}
\usepackage{graphicx}
\usepackage{caption}
\usepackage{subcaption}
\usepackage{placeins}
\usepackage{wrapfig}

\usepackage[margin=2.5cm]{geometry}

\usepackage{verbatim}

\usepackage{hyperref} 

\numberwithin{equation}{section}
\numberwithin{figure}{section}

\newtheoremstyle{thm-style-oskari}
{7pt}      
{7pt}      
{\itshape} 
{}         
{\scshape} 
{.}        
{.5em}     
{}         

\theoremstyle{thm-style-oskari}
    \newtheorem{theorem}{Theorem}[section]
    \newtheorem{proposition}[theorem]{Proposition}
    \newtheorem{corollary}[theorem]{Corollary}
    \newtheorem{lemma}[theorem]{Lemma}
    \newtheorem{definition}[theorem]{Definition}
    
    \newtheorem{convention}[theorem]{Convention}
    \newtheorem{remark}[theorem]{Remark}

\newenvironment{Proof}[1][Proof]{\begin{proof}[\sc{#1}]}{\end{proof}}

\newcommand{\NTheorem}[2] {
        \begin{theorem}[#1] \label{thr:#1}
                #2
        \end{theorem}
        }

\newcommand{\NLemma}[2] {
        \begin{lemma}[#1] \label{lmm:#1}
                #2
        \end{lemma}
        }

\newcommand{\NCorollary}[2] {
        \begin{corollary}[#1] \label{crl:#1}
                #2
        \end{corollary}
        }

\newcommand{\NDefinition}[2] {
        \begin{definition}[#1] \label{def:#1}
                #2
        \end{definition}
        }

\newcommand{\bels}[2] {
        \begin{equation} \label{#1} \begin{split} 
                #2 
        \end{split} \end{equation}
        }
\newcommand{\bea}[1]{
	\begin{align*}
		#1
	\end{align*}
	}

\definecolor{olivegreen}{rgb}{0,0.6,0.1}

\newcommand{\Lp}[1]{\mathrm{L}^{\!{#1}}}					

\newcommand{\Ind}{\mathbbm{1}}

\newcommand{\vect}[1]{\mathbf{#1}}					

\newcommand{\1} {\mspace{1 mu}}
\newcommand{\2} {\mspace{2 mu}}
\newcommand{\msp}[1] {\mspace{#1 mu}}

\newcommand{\la} {\langle}
\newcommand{\ra} {\rangle}
\newcommand{\avg}[1] {\la #1 \ra}

\newcommand{\eps}{\varepsilon}

\newcommand{\dist} {\mrm{dist}}                   
\DeclareMathOperator*{\Spec}{Spec}						

\newcommand{\mrm}[1] {\mathrm{#1}}
\newcommand{\mcl}[1] {\mathcal{#1}}

\newcommand{\brm}[1] {\boldsymbol{\mathrm{#1}}}

\newcommand{\wti}[1] {\widetilde{#1}}
\newcommand{\wht}[1] {\widehat{#1}}
\newcommand{\ul}[1] {\underline{#1}}

\newcommand{\EE} {\mathbbm{E}}
\newcommand{\PP}  {\mathbbm{P}}

\DeclareMathOperator{\supp} {supp}

\newcommand{\ins} {\msp{1}\in\msp{1}}

\newcommand{\abs}[1]{\lvert #1 \rvert}
\newcommand{\absb}[1]{\big\lvert #1 \big\rvert}
\newcommand{\absB}[1]{\Bigl\lvert #1 \Bigr\rvert}
\newcommand{\absbb}[1]{\biggl\lvert #1 \biggr\rvert}

\newcommand{\norm}[1]{\lVert #1 \rVert}
\newcommand{\normb}[1]{\big\lVert #1 \big\rVert}

\newcommand{\floor} [1] {      \lfloor {#1}        \rfloor}

\newcommand{\ceil}  [1] {      \lceil  {#1}         \rceil}

\newcommand{\R} {\mathbb{R}}
\newcommand{\C} {{\mathbb{C}}}

\newcommand{\N} {\mathbb{N}}
\newcommand{\Z} {\mathbb{Z}}

\newcommand{\Cp} {\mathbb{H}}

\newcommand{\T}{\mathbb{T}}
\newcommand{\Tc} {\mathbb{S}}
\newcommand{\DD}{\mathbb{D}}
\newcommand{\KK}{\mathbb{K}}

\newcommand{\sett}[1] { \{ {#1} \} }

\newcommand{\setb}[1] { \bigl\{ {#1} \bigl\} }

\newcommand{\setbb}[1] { \biggl\{\, {#1} \,\biggr\} }

\newcommand{\genarg} {{\,\bullet\,}}  

\newcommand{\dif} {\mathrm{d}}
\newcommand{\cI} {\mathrm{i}}
\newcommand{\nE} {\mathrm{e}}
\newcommand{\Ord} {\mathcal{O}}

\renewcommand{\Im}{\mathrm{Im}}
\renewcommand{\Re}{\mathrm{Re}}

\newcommand{\tsfrac}[2] {{\textstyle \frac{#1}{#2}}}

\newcommand{\Sx} {\mathfrak{X}} 
\newcommand{\Px} {\pi} 

\begin{document}
\renewcommand{\thefootnote}{\fnsymbol{footnote}}

\title{\bf Local spectral statistics of Gaussian matrices with correlated entries}
\author{
\begin{minipage}{0.3\textwidth}
 \begin{center}
Oskari H. Ajanki\footnotemark[1]\\
\footnotesize 
{IST Austria}\\
{\url{oskari.ajanki@iki.fi}}
\end{center}
\end{minipage}
\begin{minipage}{0.3\textwidth}
\begin{center}
L\'aszl\'o Erd{\H o}s\footnotemark[2]  \\
\footnotesize {IST Austria}\\
{\url{lerdos@ist.ac.at}}
\end{center}
\end{minipage}
\begin{minipage}{0.3\textwidth}
 \begin{center}
Torben Kr\"uger\footnotemark[3]\\
\footnotesize 
{IST Austria}\\
{\url{torben.krueger@ist.ac.at}}
\end{center}
\end{minipage}
}

\date{\today}

\maketitle
\thispagestyle{empty} 

\footnotetext[1]{Partially supported by ERC Advanced Grant RANMAT No.\ 338804, and SFB-TR 12 Grant of the German Research Council.}
\footnotetext[2]{Partially supported by ERC Advanced Grant RANMAT No.\ 338804.}
\footnotetext[3]{Partially supported by ERC Advanced Grant RANMAT No.\ 338804, and SFB-TR 12 Grant of the German Research Council}
	
\renewcommand{\thefootnote}{\fnsymbol{footnote}}

\begin{abstract}
We prove optimal local law, bulk universality and non-trivial decay for the off-diagonal elements of the resolvent for a class of translation invariant Gaussian random matrix ensembles with correlated entries. 
\end{abstract}
\vspace{0.2cm}
{\bf AMS Subject Classification (2010):} \texttt{60B20}, \texttt{15B52}

\section{Introduction}

Most rigorous works on random matrix ensembles concern either Wigner matrices with independent entries \cite{EYBull,TV} (up to the  real symmetric or complex  hermitian symmetry constraint), or invariant ensembles where the correlation structure of the matrix elements is very specific.
Since the existing  methods to study Wigner matrices  heavily rely on independence, only very few results are available on ensembles with correlated entries, see \cite{KhorunzhyPastur94, BKV1996, Chakrabarty2014, BKV} for the Gaussian case.
The global semicircle law in the non Gaussian case  with (appropriately) weakly dependent entries has been established via moment method in \cite{Schenker2005} and via resolvent method in \cite{GNT}. A similar result for sample covariance matrices was given in \cite{ORourke2012}.
All these works establish limiting spectral density on the macroscopic scale and in models where the dependence is sufficiently weak so that the limiting density of states coincides with that of the independent case. A more general correlation structure was explored in \cite{AZdep} with a nontrivial  limit density, but still only on the global scale, see also \cite{BMP2013}.
We also mention the very recent proof of the local semicircle law and bulk universality for the adjacency matrix of the $d$-regular graphs \cite{BKY, BHKY} which has  a completely different specific  correlation (due to the requirement that every row contains the same number of ones).

In this paper we consider  self-adjoint   Gaussian random matrices $ \brm{H} $ with correlated entries.
We assume that  $ \brm{H}$ is of the form $ \brm{X} + \brm{X}^\ast $ where the elements of $ \brm{X} $ have a translation invariant correlation structure. 
Our main result is the optimal local law for $ \brm{H} $, i.e., we  show that the empirical eigenvalue measure of $ \brm{H} $ converges to a deterministic probability density $ \rho $ all the way down to the scale $ N^{-1} $, the typical distance between eigenvalues,  as the dimension $ N $ of $ \brm{H} $ increases.
We also  find that the off diagonal elements of the resolvent $ \brm{G}(z) := (\brm{H}-z)^{-1}$ with $\Im \2z >0$  in the canonical basis are not negligible (unlike in the independent case) and in fact they inherit their decay from the correlation of the matrix elements.
As a simple consequence of the local law we get bulk universality.
Furthermore, we provide sufficient conditions for the asymptotic eigenvalue density $ \rho $ to be supported on a single interval with square root growth at both ends.

The proofs rely on the key observation that the (discrete) Fourier transform $ \wht{\brm{H}} = (\1\wht{h}_{\phi\1\theta}) $ of a translation invariantly correlated  self-adjoint  random matrix $\brm{H} $  has independent entries up to an additional symmetry (cf. Lemma \ref{lmm:Fourier transform} below). Thus, our recent results \cite{AEK2} on the local law and bulk universality of Wigner type matrices with a general variance matrix can be applied. Some modifications to accommodate this extra symmetry are necessary in the proofs, but they do not influence the final result.
The upshot is that in the Fourier space the diagonal elements of $ \wht{\brm{G}}(z)$
 approximately satisfy the equation
 \bels{Geq}{
-\1\frac{1}{\wht{G}_{\phi\phi}(z)} \,\approx\, z + \sum_\theta s_{\msp{-1}\phi\1\theta}\1 \wht{G}_{\theta\1\theta}(z)
\,,  \qquad s_{\msp{-1}\phi\1\theta} := \EE\,\abs{\1\wht{h}_{\phi\1\theta}}^2 \,,
} 
which constitutes a small perturbation of the \emph{Quadratic Vector Equation (QVE)},
\bels{qve}{
-\1\frac{1}{m_\phi(z)} \,=\, z + \sum_\theta s_{\phi\1\theta}\1 m_\theta(z)
\,,
}
that was extensively analysed in \cite{AEK1}. 
Since the matrix $ \brm{S} = (s_{\phi\1\theta}) $ is typically not stochastic, the components $ m_\phi(z) $ of the solution genuinely depend on $ \phi $. We establish natural conditions on the correlation structure of $\brm{H}$ that guarantee that the recently developed theory \cite{AEK1} on QVEs is applicable.
 In particular, the stability of the QVE  implies  that the solutions of \eqref{Geq} and \eqref{qve} are close, i.e., $ \wht{G}_{\phi\phi}(z) = m_\phi(z) + o(1)$, even for  spectral parameters  $z$ very close to the real axis, down to the scale  $ \Im\,z \gg N^{-1} $. This yields the local law  for the eigenvalue density of $\wht{\brm{H}}$. 
Moreover, the anisotropic law from \cite{AEK2}, applied to $\wht{\brm{H}}$,  translates directly into a precise asymptotics for any matrix elements of the resolvent in the canonical basis: 
\[
G_{xy}(z) 
\,=\,
\frac{1}{N}\sum_{\phi,\theta} \nE^{-\1\cI\12\1\pi\1(\1\phi\1x \2-\2 \theta \1y)}\,\wht{G}_{\msp{-1}\phi\1\theta}(z)
\;\approx\;
\frac{1}{N}\sum_\phi \nE^{-\1\cI\12\pi\1\phi\2(x\1-\1y)}\, m_{\phi}(z)
.
\]
The off-diagonal decay  of the entries  of $\brm{G}(z)$ thus follows from smoothness properties of $m_\phi(z)$ in the variable $\phi$.
We show that, in turn, this smoothness follows from the decay conditions on the correlation structure of $\brm{H} $. 
Finally, we  prove bulk universality of the local spectral statistics of $\brm{H}$ by using the analogous result from \cite{AEK2} for $\wht{\brm{H}} $ and the fact that $\brm{H}$ and $\wht{\brm{H}} $ are isospectral.

Gaussian random matrices with translation invariant covariance structure have been 
analyzed earlier and it has also been realized that the equation \eqref{qve} via Fourier transform plays a key role in identifying the limiting density of eigenvalues, see Khorunzhy and Pastur \cite{KhorunzhyPastur94, PasturESIrev}, Girko \cite{Girko-book}, as well as  Anderson and Zeitouni in \cite{AZdep}. These works, however, were concerned only with the density on macroscopic scales.
The off-diagonal decay of the resolvent and the bulk universality require  much more detailed information. The current paper in combination with \cite{AEK1} and \cite{AEK2} presents such a precise analysis.

\section{Set-up and main results}
\label{sec:Set-up and main results}

Consider a real symmetric or complex hermitian random matrix,
\bels{}{
\qquad
\brm{H} = (h_{ij})_{i,j\ins\T}
\,,
}
indexed by the large discrete torus of size $ N $, 
\bels{}{
\T := \Z/N\Z
\,.
}
We assume that the matrix is centered, i.e.,
\begin{subequations} 
\label{H:first two moments in x-space}
\bels{EE h = 0}{ 
\EE\,h_{ij} = 0\,,\quad\forall\, i,j \in \T
\,,
}
and that the elements $ h_{ij} $ are jointly Gaussian. 
The covariances of the elements of $ \brm{H} $ are specified by two self-adjoint matrices $ \brm{A} =  (a_{ij})_{i,j\ins\T} $ and  $ \brm{B} =  (b_{ij})_{i,j\ins\T} $, through
\bels{EE H H = A + B}{
\EE\,h_{ij}\overline{h}_{kl}
\;=\;
\frac{1}{N}(\2a_{i-k\1,\1j-l} +\, b_{i-l\1,\1j-k})
\,,\qquad\qquad
\forall\,i,j,k,l\in \T
\,.
}
\end{subequations}
Here the subtractions in $ i-k $ and $ j-l $, etc., are done in the torus $ \T $. 
Let us also denote the graph distance of $ x \in \T $ from the special point $ 0 \in \T $ by $ \abs{x} $.
We remark that any random matrix of the form $ \brm{H} = \brm{X} + \brm{X}^\ast $, where $ \brm{X} = (x_{ij})_{i,j\in\T} $ is centred and translation invariant in the sense that $(x_{i+k,j+l})_{i,j\in\T} $ has the same law as $ \brm{X} $ for any fixed shift $(k,l) \in \T^2 $, has the correlation structure \eqref{EE H H = A + B}.

The following properties of $ \brm{A} $ are needed to prove our main results:
\begin{itemize}
\item[({\bf D1})] 
{\bf Power law decay:} There is a positive integer  $ \kappa  $, such that
\bels{beta-summability of a_xy}{
\sum_{x,y\ins \T} (\11+\abs{x}+\abs{y}\1)^{\1\kappa} \abs{\1a_{xy}} \;\leq\; 1
\,.
} 
\item[({\bf D2})] 
{\bf Exponential decay:} There is a constant $ \nu > 0 $ such that
\bels{exponential decay of A}{
\abs{\1a_{xy}} \,\leq\, \nE^{-\1\nu\2(\2\abs{x}\1+\2\abs{y}\1)}
\,,\qquad
\forall\,x,y \in \T
\,.
} 
\item[({\bf R1})]
{\bf Non-resonance:} There is a constant $ \xi_1 >0 $, such that
\bels{non-resonance condition for A}{
\sum_{x\ins\T} \nE^{\1\cI\12\pi\1\phi\1x}a_{x\10} \;\ge\; \xi_1 
\,,\qquad
\forall\,\phi \in [\10,1\1]
\,.
}
\item[({\bf R2})] 
{\bf Strong non-resonance:} There is a constant $ \xi_2 >0 $, such that
\bels{strong non-resonance condition for A}{
\sum_{x,\1y\ins\T} \nE^{\1\cI\12\pi\1(\1x\1\phi\2-\2y\1\theta\1)}a_{xy} \;\ge\; \xi_2 
\,,\qquad
\forall\,\phi,\theta \in [\10,1\1]
\,.
}
\end{itemize}
In general the solution of the QVE \eqref{qve} specifying the asymptotic density of the states for $ \brm{H} $ may be neither bounded nor stable (cf. Section 9 of \cite{AEK1}). We will show that certain combinations of the above conditions  exclude these issues.

The restrictions on the correlation structure are quantified by the $ N$-independent {\bf model parameters} $ \nu,\kappa,\xi_1,\xi_2 $ appearing above.
We remark that the normalization of \eqref{beta-summability of a_xy} and \eqref{exponential decay of A} is chosen for convenience, e.g., we could replace $ 1 $ on the right hand side of \eqref{beta-summability of a_xy} by some finite constant.
The set of model parameters depends on our assumptions, e.g., if only ({\bf D1}) and ({\bf R2}) are assumed, then  $ \nu$ and $ \xi_2 $ are the model parameters. 
We allow constants appearing in the statements to depend on the model parameters. 

For compact statements of our results we define the notion of stochastic domination, introduced in  \cite{EKY} and \cite{EKYY}. This notion is designed to compare sequences of random variables in the large $N$ limit up to small powers of $N$ on high probability sets. 

\NDefinition{Stochastic domination}{
Suppose $N_0: (0,\infty)^2\to \N$ is a given function, depending only on the model parameters, as well as on an additional {\bf tolerance exponent} $ \gamma \in (0,1)$. For two sequences, $\varphi=(\varphi^{(N)})_N$ and $\psi=(\psi^{(N)})_N$, of non-negative random variables we say that $\varphi$ is {\bf stochastically dominated} by $\psi$ if for all $\eps>0$ and $D>0$,
\bels{}{
\PP\Bigl( \1\varphi^{(N)} \1>\1 N^\eps \1 \psi^{(N)}\Bigr)\;\leq\; N^{-D},\qquad N\2\geq\2 N_0(\eps,D)\,.
}
In this case we write $\varphi \prec \psi$. 
}

Let us denote the upper complex half plane and the discrete dual torus of $ \T  $ by 
\[
\Cp := \setb{z \in \C : \Im\,z > 0}\,, 
\qquad\text{and}\qquad
\Tc := N^{-1}\T\,,
\]
respectively.
It was shown in \cite{AEK1} that the Quadratic Vector Equation ({\bf QVE})
\bels{QVE in Tc}{
-\frac{1}{m_\phi(z)} \,=\, z + \sum_{\theta \ins \Tc}\wht{a}_{\phi\1\theta}\2m_\theta(z)
\,,
}
where 
\bels{def of wht-a}{
\wht{a}_{\phi\1\theta} := \frac{1}{N}\!\sum_{x,\1y\ins\T} \!\nE^{\1\cI\12\pi\1(\1x\1\phi\2-\2y\1\theta\1)}\,a_{xy}
\,,
}
has a unique solution $ \brm{m}(z) = (m_\phi(z))_{\phi\in\Tc} $ in $ \Cp^{\1\Tc} $, for every $ z \in \Cp $. 

Our main result is the optimal local law and the decay estimate for the off-diagonal resolvent entries. 
These are stated in terms of the resolvent $ \brm{G}(z)  = (G_{ij}(z))_{i,j\ins\T} $, 
\[
\brm{G}(z) \,:=\, (\1\brm{H}-z)^{-1} 
\,.
\]

\begin{theorem}[Local law for Gaussian matrices with correlated entries]
\label{thr:Local law for Gaussian matrices with correlated entries}
Suppose $ \brm{A} $ is either exponentially decaying \emph{({\bf D2})} and non-resonant \emph{({\bf R1})}, or decays like a power law \emph{({\bf D1})} and is strongly non-resonant \emph{({\bf R2})}. 
Then for any tolerance exponent $ \gamma \in (0,1) $ and uniformly for all $ z \in \R + \cI\2[\1N^{\gamma-1},\2\infty\1) $
\begin{subequations}
\begin{align}
\label{exponential off-diagonal decay of G}
\max_{x,\1y \1\in\2 \T}\,\absB{\1G_{xy}(z)-q_{x-y}(z)\1} \,
&\prec\, 
\sqrt{\frac{\Im\,q_0(z)}{N\,\Im\,z}} + \frac{1}{N\,\Im\,z}
\\
\label{Dep Gauss: local law for DS}
\absbb{\frac{1}{N}\mrm{Tr}\,\brm{G}(z) \2-\, q_0(z)}
\,&\prec\,\frac{1}{N\, \Im\,z}
\,,
\end{align}
\end{subequations}
where 
\bels{def of q_x(z)}{
q_x(z) := \frac{1}{N}\sum_{\phi\ins\Tc}\nE^{-\cI\12\pi\1x\1\phi}m_\phi(z)
\,,
\qquad
x \in \T
\,.
}
The vector $ \brm{q}(z) = (q_x(z))_{x\in \T}$ 
inherits the decay type (exponential vs. power law) from $\brm{A} $, 
in the sense that 
\bels{decay of q_x(z)}{ 
\qquad \abs{\1q_x(z)} 
\;\leq\; 
C\begin{cases}
\displaystyle
\,\abs{x}^{-\1\kappa} +\, N^{-1/2}
\quad &
\text{when \eqref{beta-summability of a_xy} holds}
\vspace{0.1cm}
\\
\displaystyle
\nE^{-\1\nu'\1\abs{\1x\1}\1}  \!+ N^{-1/2}&\text{when \eqref{exponential decay of A} holds}
\end{cases}
\qquad \forall\1x \in \T
\,,
}  
with the constants $ C > 0 $ and $ \nu' >0$ depending only on the model parameters. 
\end{theorem}
 
Generally the off-diagonal resolvent entries are not negligible even though \eqref{decay of q_x(z)} states only an upper bound.
In many cases matching lower bounds can be obtained. For example,
for the special model with correlation
$ a_{xy} := \nE^{-\nu\2(\abs{x}+\abs{y})} $
the QVE reduces to a simple scalar equation
since $a_{xy}$ factorizes. An elementary calculation shows that
in this case
as $N\to\infty$,
\[
q_x(z)\,\to\, Q(z) \1\lambda(z)^{\1\abs{x}}
\,,
\qquad \abs{x} \ge 1
\,,
\]
for some $ \lambda(z), Q(z) \in \C $ with $ 0 < \abs{\lambda(z)} < 1 $.

Note that in the general setting of Theorem 2.2 the function $ \pi^{-1}\Im\,q_0(z)  $ is the harmonic extension of the even probability density
\bels{def of density of states}{
\rho(\tau) \,:=\, \lim_{\eta \downarrow 0} \frac{1}{\pi N}\sum_{\phi \ins \Tc} \Im\,m_\phi(\tau+\cI\1\eta)
\,, 
\qquad \tau \in \R
\,,
}
to the upper half plane.
From \eqref{Dep Gauss: local law for DS} it follows that the empirical spectral measure of $ \brm{H} $ approaches the measure with the Lebesgue density $ \rho $ as $ N \to \infty $.
In fact, using a comparison argument (cf. Theorem 1 of \cite{BMP2013}) this global convergence result extends also to non-Gaussian translation invariant random matrices satisfying \eqref{H:first two moments in x-space}.
By applying the general theory for QVEs from \cite{AEK1} we are able to say more about the function $ \brm{q} : \Cp \to \C^{\T} $, and the associated even probability density $ \rho : \R \to [\10\1,\infty) $,

\begin{proposition}[Regularity of $ \rho $ and $ q_x $]
\label{prp:Regularity of q_x}
If $\brm{A} $ satisfies either \emph{({\bf D1})} and \emph{({\bf R2})}, or \emph{({\bf D2})} and \emph{({\bf R1})}, then there exists $ \beta \sim 1 $ and three constants $ C_0,c_1,C_2 > 0 $, depending only on the model parameters, such that   
$ \supp\,\rho = [-\beta,\2\beta\1] $, and  
\bels{edge behaviour of rho}{
\rho(-\beta+\omega\1) 
\,=\, 
\rho(\2\beta-\omega\1)
\,=\,
C_0\2\omega^{\11/2} \!+\2 \epsilon(\omega)
\,,\qquad
\omega \ge 0
\,,
}
where $ \abs{\1\epsilon(\omega)} \leq C_2\2\omega $.
Moreover, for an arbitrary $ \delta > 0 $, $ \rho(\tau) \,\ge\, c_1\2\delta^{1/2} $ whenever $ \abs{\tau} \leq \beta-\delta\, $. 
The function $ \brm{q} : \Cp \to \C^\T $ is analytic and it can be analytically extended to $ \R \backslash \sett{\pm\1\beta} $. In particular, the density $\rho$ is real analytic away from $\pm \beta$, the
edges of its support.
\end{proposition}

We remark that there are no explicit conditions on the correlation matrix $ \brm{B} $ in either Theorem \ref{thr:Local law for Gaussian matrices with correlated entries} or Proposition \ref{prp:Regularity of q_x}.
However, $ \brm{A} $ and $ \brm{B} $ are related. For example, if $ \brm{H} $ is real valued then $ \brm{A} = \brm{B} $. The Fourier transforms of $ \brm{A} $ and $\brm{B} $ must satisfy certain compatibility relations (cf. the proof of Corollary \ref{crl:Bulk universality}) which are equivalent to positive definiteness of the corresponding covariance matrices.

Similarly, as in the case of Wigner type matrices the local law implies the bulk universality for Gaussian matrices with correlated entries. However, the $ q $-fullness condition (Definition 1.14 in \cite{AEK2}) is replaced by a different non-generacy condition.

\begin{corollary}[Bulk universality]
\label{crl:Bulk universality}
Assume $ \brm{A} $ satisfies \emph{({\bf D1})} and  either of the following holds:
\begin{itemize}
\item 
$ \brm{H} $ is {\bf real symmetric} and $ \brm{A} $ is strongly non-resonant \emph{({\bf R2})}; 
\item 
$ \brm{H} $ is {\bf complex hermitian}, and there is a constant $ \xi_3 > 0 $ such that  
\bels{Hermitian case: A vs B in p-coords}{
\absb{\2\wht{b}_{\phi\1\theta}}^2 
\,<\, 
\Bigl(\2\wht{a}_{\phi,\theta}-\frac{\xi_3\!}{N}\2\Bigr)_+\,\Bigl(\2\wht{a}_{-\phi,-\theta}-\frac{\xi_3\!}{N}\2\Bigr)_+\,, 
\qquad\forall\1\phi,\theta \in [\10,1\1]\,, 
}
where $ \wht{b}_{\phi\1\theta} $ is defined analogously to $ \wht{a}_{\phi\1\theta} $ in \eqref{def of wht-a}, and $ \tau_+ := \max\sett{0,\tau} $, for $ \tau \in \R $.
\end{itemize}
Then for any parameter $\rho_0 > 0 $ and a smooth compactly supported function $ F:\R^n \to \R$, $ n \in \N $, there exist constants $ c,C> 0 $, depending only on $ \rho_0,\kappa $, the function $ F $, and either $ \xi_2 $ or $ \xi_3 $, such that for any $\tau \in \R $ with $\rho(\tau)\geq \rho_0 $ the local eigenvalue distribution is universal, 
\[
\absbb{\,\EE\2F\Bigr( \bigl(\2N\rho(\lambda_{i(\tau)}\msp{-1})(\lambda_{i(\tau)}-\lambda_{i(\tau)+j})\2\bigr)_{j=1}^n \Bigr)
-\,
\EE_{\rm G}\2F\Bigl( \bigl(\2N\rho_{\rm sc}(0)(\lambda_{\ceil{N/2}}-\lambda_{\ceil{N/2}+j})\2\bigr)_{j=1}^n \Bigr)
\,}
\,\leq\,
C\2N^{-c}.
\]
Here, $\EE_{\rm G}$ denotes the expectation with respect to the standard Gaussian ensemble, i.e., with respect to GUE and GOE in the cases of complex hermitian and real symmetric $\vect{H}$, respectively, and $\rho_{\rm sc}(0)=1/(2\1\pi)$ is the value of Wigner's semicircle law at the origin. 
\end{corollary}

Let us introduce the notations $ \norm{\brm{v}}_\infty := \max_i \abs{v_i} $ and $ \brm{v} \cdot \brm{w} = \sum_i \overline{\2v_i}\2w_i $ for $ \brm{v},\brm{w} \in \C^\T $.

\NCorollary{Delocalization of eigenvectors}{
Let $ \brm{u}^{(i)} \in \C^N $ be the normalized eigenvector of $\vect{H}$ corresponding to the eigenvalue $\lambda_i$. 
All eigenvectors are delocalized in the sense that for any deterministic unit vector $\brm{b}\in \C^N$ we have
\[ 
\absb{\2\brm{b}\cdot \brm{u}^{(i)}} \,\prec\, N^{-1/2}
\,.
\]
In particular, the eigenvectors are completely delocalized, i.e., $\norm{\1\vect{u}^{(i)}}_\infty\prec N^{-1/2}$.
}

The following result shows a practical way to construct real symmetric random matrices with translation invariant correlation structure. A similar, but slightly more complicated convolution representation exists for complex hermitian random matrices.

\NLemma{Linear filtering}{
Suppose a real symmetric matrix $ \brm{A} $ satisfies the Bochner type condition
\bels{Bochner in x-space}{
\sum_{i,\1j,\1k,\1l\ins\T} \!\overline{w}_{ij}\2a_{i-k\1,\2j-l}\,w_{kl} \;\ge\; 0
\,, 
} 
for arbitrary matrices $ \brm{W} = (w_{ij})_{i,j\in\T} $.
Then the random matrix $ \brm{H} $ defined as the convolution,
\bels{h_ij as convolution}{
h_{ij} := \sum_{k,l\in\T}r_{i-k,j-l}\,v_{kl}
\,,
}
of a GOE random matrix $ \brm{V} = (v_{ij})_{i,j\in\T}$, and the filter matrix $ \brm{R} = (r_{ij})_{i,j\in\T}$, defined by
\bels{def of filter R}{
r_{xy} 
\,:=\,
\frac{1}{N^{1/2}\!}\sum_{\phi,\theta \in \Tc} \nE^{-\cI\12\pi\1(x\1\phi \2-\2 y\1\theta)}\sqrt{\,\wht{a}_{\phi\1\theta}} 
\;,
}
has the correlation structure \eqref{H:first two moments in x-space} with $ \brm{B} = \brm{A} $.
}
This lemma is proven at the end of Subsection \ref{sec:Proofs for local law and bulk universality}.
We introduce the following conventions and notations used throughout this paper.

\begin{convention}[Constants and comparison relation]
\label{conv:Constants and comparison relation}
Symbols  $c, c_1, c_2,\dots $ and $ C, C_1, C_2, \dots $ denote generic positive and finite constants that depend only on the model parameters. They have a local meaning within a specific proof.  
For two arbitrary non-negative functions $ \varphi $ and $ \psi $ defined on some domain $ U $, we write
$ \varphi \,\lesssim\, \psi $, or equivalently $ \psi \,\gtrsim\, \varphi $, if $ \varphi(u) \leq C\2\psi(u) $, holds for all $ u \in U $.  
The notation $ \psi \sim \varphi$ is equivalent to both $\psi \lesssim \varphi$ and $\psi \gtrsim \varphi$ holding at the same time. 
In this case we say that $\psi$ and $\varphi$ are {\bf comparable}. In general the relation $ \gtrsim $ is called  the {\bf comparison relation}.
We also write $\psi= \varphi +\Ord(\vartheta)$ if $\abs{\psi-\varphi} \lesssim \vartheta$.
\end{convention}

\subsection{Structure of the proof}

The proof of Theorem \ref{thr:Local law for Gaussian matrices with correlated entries} splits into three separate parts. In the first part we show how to make $ \brm{H} $ into an almost Wigner type matrix by changing basis. 
In the second part we describe how the proofs for Wigner type matrices in \cite{AEK2} are modified in order to accommodate some extra dependence in the transformed matrix.
In the third part we show that the assumptions on the correlation matrix $ \brm{A} $ imply that the QVE \eqref{QVE in Tc} has a bounded and sufficiently regular solution $ \brm{m} $ using the general theory developed in \cite{AEK1}. 
Finally, in the last section we combine the results of the three steps and prove Theorem \ref{thr:Local law for Gaussian matrices with correlated entries}.

\section{Mapping $\brm{H} $ into Wigner type matrix by change of basis}

The (discrete) Fourier transforms of a matrix $ \brm{T} = (t_{ij})_{i,j\in\T} $ is another matrix $ \wht{\brm{T}} = (\1\wht{t}_{\phi\1\theta})_{\phi,\theta \in\Tc} $ defined by
\bels{def of F-transform}{
\wht{t}_{\phi\1\theta} \,:=\, \frac{1}{N}\msp{-6}\sum_{\,x,\1y\ins\T}\msp{-6}\nE^{\cI\12\pi\1(\phi\1x-\theta\1y)}\2t_{xy}
\,.
}   
Since the mapping $ \brm{T} \mapsto \wht{\brm{T}} $ corresponds to the conjugation by the unitary matrix $ \brm{F} = (f_{\phi,y})_{\phi\ins\Tc,y\ins\T} $, with elements
\bels{}{
f_{\phi\1y} :=
N^{-1/2}\nE^{\1\cI\12\pi\1\phi\1y}
\,,
\qquad
\phi \in \Tc\,,\;x \in \T
\,, 
}
the matrices $ \brm{T} $ and $ \wht{\brm{T}} = \brm{F}\brm{T}\brm{F}^\ast $ have the same spectrum:
\[
\Spec(\1\brm{T}\1) \,=\, \Spec(\1\wht{\brm{T}}\1)
\,.
\]

In the following we analyze random matrices which have independent entries modulo two reflection symmetries.

\begin{definition}[$ 4$-fold correlated ensemble]
\label{def:4-fold correlated ensemble}
A random matrix $ \brm{H} $ indexed by a torus is {\bf 4-fold correlated} if $ h_{ij} $ and $ h_{kl} $ are independent unless 
\bels{kl sim ij}{ 
(\1k,l\1) \,\in\, \setb{(\1i,j\1),(\2j,i\1),(-i,-j),(-j,-i)} 
\,.
}
\end{definition}

The next result shows that the discrete Fourier transform maps Gaussian translation invariant random matrices into Wigner type random matrices with an extra dependence. This connection was first realized by Girko \cite{Girko-book} and Khorunzhy and Pastur \cite{KhorunzhyPastur94}. It has been later used in \cite{AZdep,BMP2013,Chakrabarty2014}. 
\begin{lemma}[Fourier transform]
\label{lmm:Fourier transform}
Let $ \brm{H} $ be a (not necessarily Gaussian) random matrix satisfying \eqref{H:first two moments in x-space}.  
Then the elements of its Fourier transform $ \wht{\brm{H}} $ satisfy
\begin{subequations}
\label{H:first two moments in p-space}
\begin{align}
\EE\,\wht{h}_{\phi\1\theta} \,&=\,0
\\
\label{H: covariance in p-space}
\EE\,\wht{h}_{\phi\2\theta}\msp{-2}\overline{\,\wht{h}_{\phi'\theta'}\msp{-15}} 
\msp{10}\;&=\;
\wht{a}_{\phi\1\theta}\,\delta_{\phi\1\phi'}\,\delta_{\theta\1\theta'} +\1 \wht{b}_{\phi\1\theta'}\,\delta_{\phi,-\theta'}\,\delta_{\theta,-\phi'}
\,,
\end{align} 
\end{subequations}
for every $ \phi,\phi',\theta,\theta' \in \Tc $.
If additionally $ \brm{H}  $ is Gaussian, then $ \wht{\brm{H}} $ is $ 4$-fold correlated.
\end{lemma}
We  remark that if $ a_{xy} $ satisfies the decay estimate \eqref{beta-summability of a_xy}, then $ \wht{a}_{\phi\1\theta},\abs{\1\wht{b}_{\phi\1\theta}} \lesssim N^{-1} \!$.

\begin{Proof}
The proof of \eqref{H:first two moments in p-space} is a straightforward computation. We omit further details. 
From \eqref{H: covariance in p-space} we see that covariances between $ \Re\,\wht{h}_{\phi\1\theta} $, $ \Im\,\wht{h}_{\phi\1\theta} $ and  $ \Re\,\wht{h}_{\phi'\theta'} $, $ \Im\,\wht{h}_{\phi'\theta'} $ can be non-zero if and only if the condition equivalent to \eqref{kl sim ij} holds. 
Since the covariance matrix captures completely the dependence between the components of a Gaussian random vector 
the statement about the independence follows trivially.
\end{Proof}

\subsection{Local law for $ 4 $-fold correlation}
\label{ssec:Local law for 4-fold correlation}

In this subsection we sketch how to prove a local law for the elements of the Fourier-transformed resolvent
\[
\wht{\brm{G}}(z) \,=\, (\1\wht{\brm{H}}-z\1)^{-1}
\,,
\]
by slightly modifying the proof for the Wigner type matrices in \cite{AEK2}. Indeed, the analysis is the same as before, but due to the extra correlation between $ (\phi,\theta\1) $ and $ (-\phi,-\theta) $ we have to remove both the rows and columns corresponding to indices $ \phi $ and $ -\phi $ from $ \wht{\brm{H}} $ in order to make it independent of a given row $ \phi $. We state a local law for a general self-adjoint random matrix with independent entries apart from a possible correlation of the entries with indices 
$(i,j)$ and $(-i,-j)$.

\begin{theorem}[Local law for $ 4 $-fold correlation]
\label{thr:Local law for 4-fold correlation}
Suppose $\brm{H}=(h_{ij})_{i,j \in \T}$ is  four-fold correlated. 
If $ \brm{H}$ fulfills the conditions of Theorem 1.6 from \cite{AEK2} and has an additional symmetry
\bels{extra symmetry for 4-fold correlated matrices}{
\EE\2 h_{i\1j} h_{-j,-i}
\,=\,0\,,\qquad\qquad i\neq j
\,,
}
then the conclusions of Theorem 1.6 from \cite{AEK2} hold.

In particular, suppose the solution of 
\bels{generic QVE}{
-\frac{1}{m_i} = z + (\brm{S}\brm{m})_i
\,,\qquad i \in \T\,,\quad z \in \Cp\,, 
}
with $ s_{ij} := \EE\,\abs{h_{ij}}^2 $, is uniformly bounded in $ i$ and $z $, and that there exists a constant $ \eps_\ast > 0$ such that for every $ \eps \in (\10\1,\eps_\ast) $ the set $  \sett{\tau \in \R: \rho(\tau) > \eps} $ is an interval. Here the density $ \rho(\tau) $ is obtained by extending 
\[ 
\rho(z) :=
\frac{1}{\pi\1N}\sum_i\Im\,m_i(z) 
\,,
\]
to the real axis.
Then for any $ \gamma > 0 $ the local law holds uniformly for every $ z = \tau + \cI\1\eta $, with $ \eta \ge N^{\gamma-1}$, and non-random  $ \brm{w}\in \C^{\T}$ satisfying $ \max_i \abs{w_i} \leq 1 $:
\begin{align}
\label{4-fold: simple local law}
\max_{i\1,\2j}\absb{\1G_{ij}(z)-m_i(z)\1\delta_{ij}\1} 
\,&\prec\,
\sqrt{\frac{\rho(z)}{\!N\1\eta}}+\frac{1}{N\1\eta}
\\
\label{4-fold: simple averaged local law}
\absbb{\frac{1}{N}\sum_{i=1}^N w_i\,\bigl(G_{ii}(z)-m_i(z)\bigr)}
\,&\prec\, 
\frac{1}{N\eta}
\end{align}
The stochastic domination depends only on $\gamma $ and $ \eps_\ast $, and the constants $ \ul{\mu} =(\mu_k), P,L,p $ appearing in the estimates contained in the assumptions of  Theorem 1.6 from \cite{AEK2}: $ \EE\,\abs{h_{ij}}^k \leq (s_{ij})^{k/2}\mu_k $, $ s_{ij} \leq 1/N$, $ (\brm{S}^L)_{ij} \ge p/N $, and $ \abs{m_i(z)} \leq P $, for all $i,j \in \T $.
\end{theorem}

The extra symmetry condition \eqref{extra symmetry for 4-fold correlated matrices} is automatically satisfied by random matrices with the covariance structure \eqref{H: covariance in p-space}, but it  is generally not needed for the local law to hold (cf. \cite{JAlt} when $ \brm{S} $ is stochastic).

\begin{Proof}[Proof of Theorem \ref{thr:Local law for 4-fold correlation}] 
We modify slightly the proof of Theorem 1.6 in \cite{AEK2}.
The independence of the entries $ h_{ij} $ and $ h_{-i,-j} $ was used  to estimate the off-diagonal resolvent entries and the perturbation $ \brm{d} = \brm{d}(z) $ of the perturbed QVE satisfied by the diagonal resolvent elements $ g_k = g_k(z) = G_{kk}(z) $,
\bels{generic perturbed QVE}{
-\2\frac{1}{g_k} \,&=\, z + (\brm{S}\brm{g})_k + d_k
\,,\qquad k \in \T\,,\quad z \in \Cp\,,
}
only in the proofs of Lemma 2.1 and Theorem 3.5 in \cite{AEK2}.

In order to generalize Lemma 2.1 of \cite{AEK2} we apply the general resolvent identity (2.9) from \cite{AEK2} to replace the entries of $ \brm{G}^{(k)} $ by the corresponding entries of $ \brm{G}^{(k,-k)}$ in the defining formula (2.2) of $ d_k $ in \cite{AEK2}.
This way we obtain a representation for $ d_k $ as a sum of terms each of which can be individually shown to be small by using either trivial bounds, or by using the large deviation estimates (2.7) similarly as in the proof of Lemma 2.1 in \cite{AEK2}.
We will not present these estimates here, since a very similar analysis was carried out in Section 5 of \cite{JAlt}.
The details for obtaining this representation for $ d_k $ in the $ 4$-fold correlated random matrices are provided in Subsection 5.1 of \cite{JAlt}. We note that \eqref{generic perturbed QVE} is equivalent to formula (5.4) in \cite{JAlt} with the symbol $ \Upsilon_k $ denoting $ d_k $.  
The off-diagonal resolvent elements are treated similarly by decoupling the dependence between specific rows of $ \brm{H} $ and the entries of $ \brm{G} $. The treatment of the reflected off-diagonal elements $ G_{i,-i} $ is simpler than in Subsection 5.2 of \cite{JAlt} since the extra symmetry \eqref{extra symmetry for 4-fold correlated matrices} makes many error terms disappear.
Instead of \eqref{extra symmetry for 4-fold correlated matrices} another symmetry,  $ h_{-x,-a} = h_{ax} $, was assumed in \cite{JAlt}. 
Since all the factors of the form $ \EE\,h_{xa}^2 $ in the error terms $ \mcl{E}^{(k)}_x $ in \cite{JAlt} first appeared in the form $ \EE\,h_{ax}h_{-x,-a} $, which is zero in our case by \eqref{extra symmetry for 4-fold correlated matrices}, 
when following the proof in \cite{JAlt}, we may replace the terms $ \EE\,h_{xa}^2 $ with zeros.

The fluctuation averaging (Theorem 3.5 of \cite{AEK2}) is extended to $ 4$-fold correlated matrices also  by slightly modifying the original proof of Theorem 4.7 in \cite{EKYY}. In particular, the arguments do not rely on the stochasticity of $ \brm{S} $ as explained in the proof of Theorem 3.5 in \cite{AEK2}.
In order to handle the extra dependence one needs to make a simple modification: The equivalence relation given within the proof needs to be generalized such that for a given index-tuple $\vect{k}=(k_1,\dots,k_p)\in \T^p$, we define $ r \sim s $ to mean that either $ k_r=k_s$ or $k_r=-k_s$. 
This means that for each 'lone index' $k$ one removes the index $-k$ in addition to $k$ from the other resolvent elements within the same monomial. 
For a more detailed description of the necessary modifications see the proof of Theorem 4.6 in \cite{JAlt}. 
\end{Proof}

\subsection{Anisotropic local law for $ 4 $-fold correlation}

In order to translate the statements of the local law in Fourier coordinates back to the original coordinates we need an \emph{anisotropic local law}.
Here we consider $ \abs{z} $ to be bounded to get simpler estimates. This condition can be easily dropped out if needed.

\NTheorem{Anisotropic law}{
Suppose $\brm{H}=(h_{ij})_{i,j \in \T}$ is a self-adjoint $ 4$-fold correlated random matrix with centered entries satisfying the local law at some fixed $ z $, satisfying $ \abs{z} \leq 10 $,
\bels{generic local law}{
\max_{i,j}\absb{\1G_{ij}-\delta_{ij}\1m_i\1} \,\prec\, \Phi
\,,
}
where the non-random constant $ \Phi $ satisfies $ N^{-1/2} \leq \Phi \leq N^\kappa $, for some $ \kappa > 0 $.
Then uniformly for all $z= \tau +\cI\1\eta \in \Cp$ satisfying $\eta \geq N^{\gamma-1}$, and all non-random unit vectors  $\brm{u},\brm{v} \in \C^{\T}$:
\bels{isoeq}{
\absb{\1\brm{u}\cdot (\1\brm{G}-\mrm{diag}(\brm{m})\1)\1\brm{v}\1}
\,\prec\,
\Phi
\,.
}
}

\begin{Proof}
The proof is a straightforward generalization of the method applied in \cite{EKYYiso} to prove anisotropic local law for random covariance matrices and general Wigner matrices. 
The proof boils down to showing that for every $ p \in 2\1\N $ there exists a constant $ C(p) $ independent of $ N $ such that for every $ \norm{\brm{v}}_{\ell^{\12}} \leq 1 $ the moment bound
\bels{anisotropic goal}{
\EE \;\absB{\sum_{\,a\neq\1 b}\!\overline{v}_a G_{\msp{-2}ab} v_b}^{\2p} 
\,\leq\,C(p)\1\Phi^{\2p}\,,
} 
holds. In the following we will denote generic constants depending only on $ p $  by $ C(p) $.
Only two minor modifications to the method used in Section 7 of \cite{EKYYiso} are needed. First, since $\brm{S}$ is not stochastic one needs to take into account that $ G_{ii}(z) $ is close to $ m_i(z) $ instead of an $ i $-independent function such as the Stieltjes transform of the semicircle law, $ m_{\mrm{sc}}(z) $. This generalisation was handled in \cite{AEK2} (cf. Theorem 1.12) where the anisotropic local law was proven for Wigner type matrices.
As the second modification we need to incorporate the extra dependencies between the matrix elements $ h_{ak}$ and $ h_{-a,-k} $ into the analysis of \cite{EKYYiso}.
To this end we walk through the key points of the arguments in \cite{EKYYiso} and point out along the way how the steps are modified to incorporate this extra dependence.

The starting point of the argument is to write the right hand side of \eqref{anisotropic goal} in the form:
\bels{EE Z^p opened up}{
\msp{-30}\EE \,\absbb{\sum_{\,b_1 \neq\1 b_2} \msp{-8}\overline{v}_{b_1}G_{b_1b_2} v_{b_2}}^p
\,=
\sum_{b_{11}\neq  b_{12}} \msp{-7}\cdots\msp{-7} \sum_{b_{p1}\ne b_{p2}}\msp{-7}
\overline{v}_{b_{11}}v_{b_{12}}\cdots\, \overline{v}_{b_{p1}}v_{b_{p2}}
\;\,
\EE\;
\prod_{k\1=\11}^{p/2} \!G_{b_{k1}b_{k2}}\msp{-14}\prod_{l=p/2+1}^p \msp{-15}G^\ast_{b_{l1}b_{l2}}
\,,
}
for an arbitrary even integer $ p $.
Let us consider a fixed summand, so that the values of the $ \brm{v} $-indices $ b_{k1},b_{k2} $ are fixed. Here the size of the expectation is naively bounded by $ \Phi^p $. However, there are $ 2\1p$-sums over the elements of the $ \ell^{\12}$-unit vector $ \brm{v}$. Since $ \norm{\brm{v}}_{\ell^{\11}} \leq N^{1/2} $ the naive size of the right hand side of \eqref{EE Z^p opened up} is hence $ N^{p/2}\Phi^p $.

The key idea of the proof is to apply recursively the general resolvent identities (cf. (2.9) in \cite{AEK2}) to express the product of resolvent entries in \eqref{EE Z^p opened up} as a sum over so-called \emph{trivial leaves} (cf. Subsection 5.10 of \cite{EKYYiso}) and the sum over $ C(p) $ terms (corresponding to the \emph{non-trivial leaves} in \cite{EKYYiso}) of the form
\bels{generic good term}{
\Gamma_{\!\brm{a},\brm{b},\brm{c},\brm{d}}\,
\prod_\beta h_{c_\beta d_\beta}
\sum^{(B)}_{\brm{i}}\sum_{\brm{j}}^{(B)}\delta(\brm{i},\brm{j})
\prod_\alpha h_{a_\alpha i_\alpha}G^{(B),\#}_{i_\alpha j_\alpha} h_{j_\alpha b_\alpha}
\,.
}
Here $ B = \sett{\pm\1b_{k 1}} \cup \sett{\pm\1b_{k\12}} $ is the set of all rows of $ \brm{H} $ that may depend on the rows indexed by the $ \brm{v} $-indices,  $ \sett{a_\alpha} \cup \sett{b_\alpha}\cup \sett{c_\beta}\cup \sett{d_\beta} \subset B $ and $ \brm{i} = (i_\alpha) $ and $ \brm{j} = (j_\alpha) $ are summed over $ \T \backslash B $ with a non-random indicator function $ \delta(\brm{i},\brm{j})$ possibly further restricting these sums. The superscripts $ \# $ on resolvent entries indicate possible hermitian conjugations.
The products in \eqref{generic good term} contain at most $ C(p) $ factors, while the symbol $ \Gamma_{\!\brm{a},\brm{b},\brm{c},\brm{d}} $ denotes a non-random function of size $ \Ord(1) $ that may depend on $ a_\alpha,b_\alpha,c_\beta,d_\beta $. 
We remark that terms of the form \eqref{generic good term} are coded by expressions $ \mcl{A}_{\brm{a}_b,\brm{a}_b}(\Gamma_{\!\zeta}) $ in the formula (5.45) of \cite{EKYYiso}.

The trivial leaves, exactly as in \cite{EKYYiso}, correspond to products of resolvent entries that remain smaller than $ \Phi^p $ even after summing over the $ \brm{v} $-indices simply  because they contain products of so many off-diagonal resolvent elements that these together compensate the factor $ N^{p/2} $ originating from the brutal $\ell^{\11} $-summation over $ \abs{v_b} $ (cf. Subsection 5.11 of \cite{EKYYiso}). 
The classification of the constituents of the product of resolvent entries into the trivial and the non-trivial leaves  relies on the concept of \emph{maximally expanded} resolvent entries (Subsection 5.3 in \cite{EKYYiso}).
For $ 4$-fold correlated matrices we redefine resolvent entries of the type $ G^{(B\backslash ab)}_{ab} $, with $ a,b \in B $, as being maximally expanded.
Here the set $ B $ plays the same role as the \emph{black vertices} in \cite{EKYYiso}.

From now on we concentrate on the non-trivial leaves of the type \eqref{generic good term}. 
The key property of these terms is that their expectation factorizes into an expectation over all the entries of $ \brm{H} $, and the expectation over all the entries of $ \brm{G} $.
From \eqref{EE h = 0} it follows that the expectation over the entries of $\brm{H}$ can be non-zero only if each entry of $ \brm{H} $ is paired with at least one of the four possible entries of $ \brm{H} $ it is not independent of. 
As a consequence, either each $\brm{v}$-index is paired with at least one other $ B $-index, or there are so many extra entries of $ \brm{H} $ compared to the number of independently summable indices in the products of \eqref{generic good term} that the small sizes $ \abs{h_{ai}} \prec N^{-1/2} $ compensate the presence of non-paired $\brm{v}$-indices.
In order to see that every term of the type \eqref{generic good term} indeed has these properties one uses the same graph expansion as in \cite{EKYYiso} to perform the relevant bookkeeping. 

The key insight about the combinatorics of the pairing of $ \brm{H}$-entries is that the number of ways to pair all $ C(p) $ of them in \eqref{generic good term} is bounded by a number only depending on $ p $, say by $ C(p)^{C(p)} $, but not on $ N $. Such a factor can be included in the constant $ C(p) $ on the right hand side of \eqref{anisotropic goal}, and is hence harmless.
Since  $ h_{bk} $ may be paired not only with itself but also with $h_{-b,-k}$, it is now possible that $ v_b $ gets paired with $ v_{-b} $. However, using  
\bels{CS for 4-fold correlated v's}{ 
\abs{\1v_{\1a}}\1\abs{\1v_{-a}} \,\leq\, \abs{\1v_{\1a}}^2+ \abs{\1v_{-a}}^2 
\,,
}
these terms can be reduced to $ \ell^2 $-norms of $ \brm{v} $.

Let us illustrate the modifications by considering the simplest leading order terms of the type \eqref{generic good term} when $ p = 2 $.  Considering the contribution of such terms to the right hand side of \eqref{EE Z^p opened up} yields
\[
\sum_{a\neq b}\sum_{c\neq d}\, \overline{v}_a v_b\2 \overline{v}_c v_d
\;m_a m_b\1 \overline{m}_c\overline{m}_d
\;\EE\!\sum_{i,j,k,l}^{[a,b,c,d\1]}
G^{[a,b,c,d\1]}_{ij}G^{[a,b,c,d\1]\1*}_{kl}
\,\EE\, h_{ai}\overline{h}_{bj}h_{ck}\overline{h}_{dl}
\,,
\]
where $ \brm{G}^{[T]} := \brm{G}^{(\1T \1\cup (-T))} $ for any set $ T \subset \T $.
Here the product $ m_a m_b\1 \overline{m}_c\overline{m}_d $ corresponds to the constant $ \Gamma_{\!\brm{a},\brm{b},\brm{c},\brm{d}} = \Ord(1) $ and the inner sums correspond to the $ \brm{i}$, $ \brm{j} $-sums in \eqref{generic good term}.
Without the $ 4$-fold correlations there are only two ways to pair the entries of $ \brm{H} $: (1) $ (a,i) = (c,k) $,$ (b,j)=(d,l) $, and  (2) $ (a,i) = (d,l)$, $ (b,j) = (c,k) $. 
On the other hand, under $ 4$-fold correlations it is possible to pair the entries in $ 9 $ different ways:  (1--4) $ (a,i) = \pm\1(c,k) $, $ (b,j)=\pm\1(d,l) $; (5--8) $ (a,i) = \pm\1(d,l) $, $ (b,j) = \pm\1(c,k) $; and (9) $ (a,i) = -\1(b,j) $, $ (c,k) = -\1 (d,l) $. Here, $ -\1(b,j) := (-b,-j) $, and the signs $ \pm $ can be chosen independently of each other.
The pairings possible without the $ 4$-fold correlation yield terms such as
\[
\sum_{a\neq b} \abs{v_a}^2\abs{v_b}^2 \frac{1}{N^2}\EE\sum_{i,j} \abs{G^{[a,b]}_{ij}}^2
\;\lesssim\;
\left(\EE\max_{i\neq j} \abs{G^{[a,b]}_{ij}}+ \frac{1}{N}\EE\max_i \abs{G^{[a,b]}_{ii}}^2\right)\norm{\brm{v}}_{\ell^{\12}}^4\,,
\]
which are stochastically dominated by $ C(p)\2\Phi^2$ since $ \abs{G^{(T)}_{ij}-\delta_{ij}m_i} \prec C(p)\2\Phi $ for any set $ T \subset \T $ satisfying $ \abs{T} \leq 2p $.  Here we have used $ \abs{G^{(T)}_{ij}-G_{ij}} \lesssim C(p)\2\Phi $ (cf. (2.10) from \cite{AEK2}) and the local law \eqref{generic local law}.
Under the $ 4$-fold correlations the pairing produces also terms such as
\bels{v_a v_-a match} {
\sum_{a,c} \overline{v}_a\2v_{-a}\2\overline{v}_c\2v_{-c}\frac{1}{N^2}\EE\sum_{i\1,\2k}\absb{G^{[a,c]}_{i,-i}}\absb{G^{[a,c]}_{k,-k}}
\,.
}
Here the off-diagonal resolvent elements are again stochastically dominated by $\Phi$. The sums over $ a $ and $ c $ can be bounded using \eqref{CS for 4-fold correlated v's} and $ \norm{\brm{v}}_{\ell^{\12}} \leq 1 $. Hence, also \eqref{v_a v_-a match} is stochastically dominated by $\Phi^2$.  
\end{Proof}

\section{Properties of QVE}

In this section we show that the assumptions on $ \brm{A} $ in our main theorems guarantee that the induced QVE \eqref{QVE in Tc} has a sufficiently regular uniformly bounded solution.
We show that the quantity $ q_{x-y}(z) $ describing the asymptotic value of the off-diagonal resolvent elements $ G_{x-y} $ (cf. \eqref{exponential off-diagonal decay of G}) has the correct decay properties in $ \abs{x-y} $ by using the regularity of the solution of the QVE. 

Let us define a function $ \wti{a} : [0,1]^2 \to \C $ as a continuous extension of the elements of $ N\wht{\brm{A}} $, 
\bels{def of wti-a}{
\wti{a}(\phi,\theta) := \sum_{k,l=0}^{N-1} a_{kl}\2e_k(\phi)\1e_{-l}(\theta)
\,,\qquad
\phi,\theta \in [\10\1,1]
\,,
}
where $ e_k : \R \to \C $ denotes the exponential function $ e_k(\phi) := \nE^{\2\cI\12\pi\1k\2\phi} $. 
Here we identified $ \T $ with the set of integers $ \sett{0,1,2,\dots,N-1} $. 
We remark that $ \wti{a}(\phi,\theta) \ge 0 $ for all $ \phi,\theta \in [0,1] $. This follows from the Bochner inequality \eqref{Bochner in x-space}.
Note that $ \wht{a}_{\phi\1\theta} = N^{-1}\wti{a}(\phi,\theta) $, if $ \phi,\theta \in \Tc $, with $ \Tc $ being canonically embedded in $ [0,1] $.  

We will now define a third non-resonance condition for a correlation matrix $ \brm{A} $ in terms of the induced integral operator $ \wti{A} $ acting on functions $ h: [\10\1,1] \to \C $,
\bels{def of wti-A}{
\wti{A}h(\phi) \,:=\, \int_0^1 \wti{a}(\phi,\theta)\1h(\theta)\2\dif \theta
\,.
}
\begin{itemize}
\item[({\bf R0})]
The integral operator $ \wti{A} $ is {\bf block fully indecomposable} (cf. Definition 1.7 of \cite{AEK1}), i.e., 
there exist two constants $ \xi_0 > 0 $, $ K \in \N $, a fully indecomposable  matrix $ \brm{Z} = (Z_{ij})_{i,j=1}^K $, with $Z_{ij} \in \sett{0,1} $, and a measurable partition $ \mcl{D} := \sett{D_j}_{j=1}^K $ of $ [0,1] $, such that for every $ 1 \leq i,j \leq K $ the following holds:
\bels{B2: Quantitative block FID condition}{
\abs{D_j} \,=\,\frac{1}{K}\,,
\qquad\text{and}\qquad
\wti{a}(\phi,\theta) \,\ge\,\xi_0\1 Z_{ij}\,,\quad\text{whenever}\quad(\phi,\theta) \in D_i \times D_j
\,.
}
\end{itemize}
If ({\bf R0}) is assumed we will treat the associated parameters $ \kappa, K,\xi_0 $ as model parameters. 
By definition ({\bf R2}) implies $ \wti{a}(\phi,\theta) \ge \xi_2 $ for every $ \phi,\theta $, and thus ({\bf R0}) holds with $ \xi_0 = \xi_2 $ and $ \mcl{D} = \sett{[0,1]} $. Assumption ({\bf R1}) does not imply ({\bf R0}), but ({\bf R1}) and ({\bf D2}) together do ( cf. Lemma \ref{lmm:R1 and D2 imply R0} below).

Instead of directly analyzing the discrete QVE \eqref{QVE in Tc} we will first establish the correct properties for the solution of the continuous version
\bels{continuous QVE}{
-\frac{1}{\wti{m}(z)} \,=\, z + \wti{A}\1\wti{m}(z)
\,, 
} 
of \eqref{QVE in Tc}. Afterwards we deduce these properties for the discrete version \eqref{QVE in Tc} as well. 
For the transition from the  discrete to the continuous version we need certain stability properties of the QVE that were established in \cite{AEK1}.

We recall several notations and results from \cite{AEK1}.
We will consider QVEs defined on a probability space $ (\Sx,\Px) $ with an operator $ S $ in two different setups.
When we discuss the discrete QVE \eqref{QVE in Tc}, the setup is 
\begin{subequations}
\label{QVE setups}
\bels{QVE setups: discrete}{
\Sx := \Tc
\,,\qquad
\Px := \frac{1}{N}\sum_{\phi \ins\Tc} \delta_\phi
\qquad\text{and}\qquad
S := N\wht{\brm{A}}\,,\;\text{ i.e., }\;
S_{\phi\1\theta} := N\,\wht{a}_{\phi\1\theta}
\,.
}
For the continuous QVE \eqref{continuous QVE} the setup is
\bels{QVE setups: continuous}{
\Sx := [\10\1,1\1]
\,,\qquad
\Px(\dif \phi) := \dif \phi
\qquad\text{and}\qquad
S := \wti{A}\,,
\;\text{ i.e., }\;
S_{\phi\1\theta} := \wti{a}(\phi,\theta)
\,.
}
\end{subequations}
In the following, all $ \Lp{p}$-norms and the scalar products are understood in the appropriate setups \eqref{QVE setups}.

\begin{lemma}[Bounded solution]
\label{lmm:Bounded solution}
If $ \brm{A} $ satisfies \emph{({\bf D1})} and \emph{({\bf R0})}, then the 
solution $ \wti{m}(z) : [0,1] \to \overline{\Cp} $ of the continuous QVE \eqref{continuous QVE} satisfies
\bels{uniform bounds of wti-m}{
\abs{\1\wti{m}(z;\phi)} \,+\,
\abs{\1\partial_\phi\wti{m}(z;\phi)}
\;\lesssim\;
1 
\,,
\qquad
\forall\,(z,\phi) \in \Cp \times [\10\1,1\1]
\,.
}
The unique solution $ \brm{m} $ to the discrete QVE \eqref{QVE in Tc} is close to $ \wti{m} $:
\bels{m and wti-m difference}{
\sup_{\phi\ins\Tc}\,\abs{\1m_\phi(z)-\wti{m}(z;\phi)} 
\,\lesssim\, 
N^{-1/2}
\,,
\qquad
\forall\,z \in \Cp
\,.
}
\end{lemma}
\begin{Proof}

We prove \eqref{uniform bounds of wti-m} first. 
In order to apply the general theory for QVEs we first show that the integral operator $ \wti{A} $ satisfies the conditions {\bf A1-A5.} from \cite{AEK1}. The qualitative properties {\bf A1.} and {\bf A2.} are trivially satisfied. 
For the property {\bf A5.} we show that the integral kernel of $ \wti{A}^{K-1} $ is bounded from below by a constant comparable to one. This follows from ({\bf R0}) since every element of the the $ (K-1)$-th power of the indecomposable matrix $ \brm{Z} $ is equal to or larger than one (cf. Proposition 4.3 of \cite{AEK1}).
For {\bf A4.} we need to show that the norm $ \norm{\wti{A}}_{2\to\infty} $ of $ \wti{A} $ as an operator from $ \Lp{2}[0,1] $ to $ \Lp{\infty}[0,1] $ is $ \Ord(1)$. This follows from \eqref{def of wti-a}, because
\bels{C^1-continuity of wti-a}{
\abs{\2\wti{a}(\phi,\theta)} + 
\abs{\2\partial_\phi \wti{a}(\phi,\theta)\1}
\,\leq\,
2\1\pi\! \sum_{x,\2y=0}^{N-1} (\11+\abs{x}\1)\2\abs{\1a_{xy}}
\;\lesssim\;
1
\,.
}
Finally, the normalization {\bf A3.} of \cite{AEK1} holds if we replace $ \wti{A} $ and $ \wti{m}(z;\phi) $ by $ \lambda\1\wti{A} $ and $  \lambda^{1/2}\wti{m}(\lambda^{1/2}z;\phi) $, respectively, with $ \lambda := \norm{\wti{A}}_{1\to 1}^{-1} $. From  \eqref{C^1-continuity of wti-a} it follows that $\lambda \sim 1 $. 

Next we show that $ m(z) $ is uniformly bounded for $ z \neq 0 $. Indeed, using \eqref{C^1-continuity of wti-a} we get
\[
\norm{\1\wti{a}(\phi_1,\genarg)-\wti{a}(\phi_2,\genarg)}_2 \,\leq\, C_2\1\abs{\phi_1-\phi_2} 
\,,\qquad
\forall\2\phi_1,\phi_2 \in [\10\1,1\1]
\,.
\]
From this it follows that
\[
\lim_{\eps\to 0}
\inf_{\phi_1\in\,[\10\1,1]}\int_0^1\!\frac{\dif \phi_2}{(\2\eps+\norm{\2\wti{a}(\phi_1,\genarg)-\2\wti{a}(\phi_2,\genarg)}_2)^2\!} \;=\; \infty 
\,.
\]
Since this implies the condition {\bf B1.} of \cite{AEK1} 
(i) of Theorem 4.1 in \cite{AEK1} is applicable in the setup \eqref{QVE setups: continuous}.
The theorem shows that $ \norm{m(z)}_\infty \leq C(\delta) $ for any $ \abs{z} \ge \delta $ with $ C(\delta) $ depending on $ \delta > 0 $.
The property ({\bf R0}) is equivalent to property {\bf B2.} in \cite{AEK1}. 
Hence by (ii) of Theorem 4.1 in \cite{AEK1} $ \wti{m}(z) $ is uniformly bounded in some neighborhood of $ z = 0 $. Combining this with the uniform bound away from $ z = 0 $ we get the uniform bound for $ \abs{\wti{m}(z;\phi)} $ for all $ z $ and $ \phi $.
In order to bound also the derivative $ \partial_\phi \wti{m}(z;\phi) $ we differentiate the continuous QVE \eqref{continuous QVE} and get
\bels{differential of wti-m}{
\partial_\phi\1 \wti{m}(z;\phi) \,=\, \wti{m}(z;\phi)^{\12}\!\int_0^1
\msp{-6}\wti{m}(z;\theta)\, \partial_\phi\wti{a}(\theta,\phi)\2\dif \theta
\,.
}
Using \eqref{C^1-continuity of wti-a} and the uniform boundedness of $ \wti{m} $ we  finish the proof of \eqref{uniform bounds of wti-m}.

Next we define
\bels{def of wti-rho}{
\wti{\rho}(\tau) := \lim_{\eta \downarrow 0}\frac{1}{\pi}\int_0^1 \Im\,\wti{m}(\tau +\cI\1\eta\1;\phi)\1\dif \phi
\,,
}
analogously to $ \rho $ in \eqref{def of density of states}, in the continuous setting. 
With Theorem 1.1 of  \cite{AEK1} we see that $ \wti{\rho} $ is a bounded probability density on $ \R $.
By the continuity \eqref{C^1-continuity of wti-a} we also have
\bels{single interval condition for wti-a}{
\sup_{D \1\subset\1 [\10\1,1\1]} 
\inf_{\substack{\phi_1 \ins D\\ \phi_2 \1\notin\1 D}} 
\normb{\2\wti{a}(\phi_1,\genarg)-\2\wti{a}(\phi_2,\genarg)}_1 
\;=\;
0
\,,
}
and hence Theorem 1.9 from \cite{AEK1} yields $ \supp \wti{\rho} \,= [-\wti{\beta}\1,\2\wti{\beta}\,] $, for some constant $ \wti{\beta} \sim 1 $.

Now we prove \eqref{m and wti-m difference}  by considering \eqref{QVE in Tc} as a perturbation of \eqref{continuous QVE}. Given $ \brm{m} $ we first embed $ \Tc $ into $ [\10,1] $ canonically, and define the piecewise constant functions
\bels{}{
g(z;\phi) \,&:=\, m_{N^{-1}\floor{N\phi}}(z)
\\
t(\phi,\theta) \,&:=\, N\,\wht{a}_{N^{-1}\floor{N\phi},\2N^{-1}\floor{N\theta}}
\,,
}
for $ \phi,\theta \in [0,1] $. 
Notice that $ t(\phi,\theta) = \wti{a}(\phi,\theta) $ and $ g(z;\phi) = m_\phi(z) $ when $ \phi,\theta \in \Tc $. Hence it is enough to estimate $ \norm{g-\wti{m}}_\infty $. Together with \eqref{C^1-continuity of wti-a} this implies
\bels{t- wti-a uniform}{
\abs{\2t(\phi,\theta)-\wti{a}(\phi,\theta)} \,\lesssim\, N^{-1}
\,,
\qquad\phi,\theta \in [0,1]
\,. 
}
In terms of these quantities \eqref{QVE in Tc} can be written as 
\bels{QVE imbedded}{
-\frac{1}{g} \,=\,z + Tg\,,
\qquad\text{where}\qquad
Th(\phi) \,:=\, \int_0^1 \!t(\phi,\theta)\1h(\theta)\1\dif \theta
\,.
}
We will now consider $ g $ as the solution of the perturbed continuous QVE 
\bels{g as a solution of perturbed cQVE}{
-\frac{1}{g}
\,=\,
z +
\wti{A}g + d
\,,\qquad\text{where}\quad
d := (T-\wti{A})\1g\,.
}
Using \eqref{t- wti-a uniform} we see that the perturbation $ d $ is indeed small:
\bels{T-wtiA:bound on d}{
\norm{d}_\infty \leq \norm{\1T-\wti{A}\1}_{2\to\infty}\norm{g}_2 \,\lesssim\, N^{-1}\norm{g}_2
\,.
} 
Clearly, $ \norm{T}_{2\to\infty} \sim 1 $ as well. Hence, we know from the general theory (cf. the bound (1.2) of Theorem 1.1 of \cite{AEK1}) that $ \norm{g(z)}_2 \leq 2/\abs{z} $. 
Using \eqref{t- wti-a uniform} we see that that for sufficiently large $ N $ the operator $ T $ is also block fully indecomposable with the same matrix $ \brm{Z} $ and the same partition $ \mcl{D} $ as $ \wti{A} $. Thus we get $ \norm{g(z)}_2 \lesssim 1 $ for all $ z $ by (ii) of Theorem 4.1 in \cite{AEK1}.
Combining this with \eqref{T-wtiA:bound on d} yields
\bels{unif bound on d}{ 
\norm{\1d(z)}_\infty \,\lesssim\, N^{-1}
\,.
}

Comparing \eqref{g as a solution of perturbed cQVE} and \eqref{continuous QVE} we show that \eqref{unif bound on d} implies that the corresponding solutions $g$ and $\widetilde m$ are close in the sense of \eqref{m and wti-m difference}. For this purpose we use the rough stability statement from Theorem 1.10 of \cite{AEK1} to get
\bels{GapEstimate1}{
\norm{\1g(z)-\wti{m}(z)}_\infty \Ind\bigl(\2\norm{g(z)-\wti{m}(z)}_\infty \leq \lambda_1 \bigr) 
\;\lesssim\;  
N^{-1}, 
\qquad 
\dist(z,\{\wti{\beta},-\wti{\beta}\})
\ge c_0
\,,
}
where $ c_0 \sim 1 $ and $ \lambda_1 \sim 1 $ are sufficiently small constants left unspecified until the end of the proof.

This means that we get stability as long as we stay away from the points $ \pm\1\wti{\beta} $. The necessary initial bound inside the indicator function is satisfied for large enough values of $|z|$, since 
\[
\norm{g(z)}_\infty + \norm{\widetilde m(z)}_\infty\,\lesssim\, |z|^{-1}\,, \qquad |z|\geq C_1
\,.
\]
Here $ C_1 $ is a sufficiently large constant. This bound follows from the Stieltjes transform representation of both the solution of the discrete and the continuous QVE (cf. Theorem 1.1 of \cite{AEK1}). We use continuity of $g$ and $\widetilde m$ in $z$ and \eqref{GapEstimate1} to propagate the initial bound from the regime of large values of $|z|$ to all $z \in \Cp$ with $\dist(z,\{\wti{\beta},-\wti{\beta}\})\geq c_0 $. In particular, \eqref{GapEstimate1} remains true even without the indicator function, i.e.,
\bels{BoundAwayFromSigmaStar}{
\norm{g(z)-\wti{m}(z)}_\infty  
\,\lesssim\,  
N^{-1}, \qquad 
\dist(z,\{\1\wti{\beta},-\wti{\beta}\})
\,\ge\, 
c_0
\,.
} 

It remains to show \eqref{m and wti-m difference} close to the edges by using that the instability at these two points is quadratic. The argument is a simplified version of the one used in a random setting in Section 4 of \cite{AEK2}. For the convenience of the reader we show a few details.
We restrict to the case $|z-\wti{\beta}|\leq c_0 $, close to the right edge. The left edge is treated in the same way. 
For the following analysis we use the stability result, Theorem 4.2 of \cite{AEK2}, in the continuous setup (cf. Proposition 8.1 in \cite{AEK1}). The theorem yields 
\bels{GapEstimate2}{
\norm{g-\wti{m}}_\infty \Ind\bigl(\2\norm{g-\wti{m}}_\infty \leq \lambda_2 \bigr) 
\;\lesssim\;
\Theta  \,+\, N^{-1}
\,,
}
where the quantity $ \Theta = \Theta(z) \ge 0 $ is continuous in $z$ and satisfies the cubic inequality
\bels{CubicAtEdge}{
\absb{\2\Theta^3+ \pi_2\1\Theta^2 +\pi_1\1\Theta\2}\2\Ind\bigl(\2\norm{g-\wti{m}}_\infty \leq \lambda_2 \bigr)  
\;\lesssim\; 
N^{-1}.
}
Here the constant $ \lambda_2 \sim 1 $ is independent of $ c_0 $.

Note that \eqref{CubicAtEdge} corresponds to (4.10) in \cite{AEK2} and (8.5) in \cite{AEK1}, respectively. 
Combining (4.11), (4.14b) and (4.5d) in \cite{AEK2}, the coefficients $\pi_k = \pi_k(z) $ of the cubic equation \eqref{CubicAtEdge} satisfy 
\bels{scaling of pi_k's}{
|\1\pi_1|\,\sim\,\abs{z-\widetilde{\beta}\1}^{1/2} \leq c_0^{1/2},
\qquad\text{and}\qquad
|\1\pi_2|\,\sim\, 1
\,,
}
provided $c_0 \sim 1 $ is sufficiently small. Since $\pi_1(z) \to 0 $ as $z \to \widetilde \beta$, by decreasing the size $c_0$ of the neighborhood we are working on, the value of $|\pi_1|$ can be made arbitrarily small. This, in turn, implies that the solution $\Theta$ of the cubic inequality \eqref{CubicAtEdge} is small,
\[
\Theta\2\Ind\bigl(\2\norm{g-\wti{m}}_\infty \leq \lambda_2 \bigr)
\;\lesssim\;
\abs{\1\pi_1} + N^{-1/2} .
\]
Using this we can make the right hand side of \eqref{GapEstimate2} smaller than $\lambda_2/2$, say, by decreasing the value of $c_0$. 
Thus, there is a gap in the possible values of the continuous function $z \mapsto\norm{g(z)-\wti{m}(z)}_\infty$, in the sense that $ \norm{g-\wti{m}}_\infty \notin (\1\lambda_2/2,\2\lambda_2) $.
Since on the boundary, $|z-\wti{\beta}\1|=c_0 $, the initial bound, $\norm{g-\wti{m}}_\infty \leq \lambda_2 $, holds by \eqref{BoundAwayFromSigmaStar}, it propagates to all $z$ with $\abs{z-\wti{\beta}}\leq c_0 $.
Thus, \eqref{GapEstimate2} and \eqref{CubicAtEdge} remain true without the indicator functions. 

It still remains to bound $ \Theta $ in \eqref{GapEstimate2}.
Since $|\pi_2|\sim 1$, we may absorb the cubic term in $\Theta$ in \eqref{CubicAtEdge}. We find that $\Theta$ satisfies 
\bels{QuadraticAtEdge}{
|\Theta^2 + \varpi \1\Theta| \,\lesssim\,N^{-1}\,,\qquad |\varpi|\,\sim\, |z-\widetilde \beta|^{1/2}
\,,
}
where  $ \varpi := \pi_1/(\pi_2+\Theta\1) $. 
From this it is easy to see that the bound $\Theta \leq N^{-1/2}$ can be propagated from the boundary $|z-\widetilde \beta|= c_0 $ inside the neighborhood $|z-\widetilde \beta|\leq c_0 $  of the right edge to give $ \Theta \lesssim N^{-1/2}$ everywhere. Using this in \eqref{GapEstimate2} without the indicator function proves the bound \eqref{m and wti-m difference} at the right edge. 
\end{Proof}

\begin{lemma}
\label{lmm:R1 and D2 imply R0} 
If $ \brm{A} $ satisfies \emph{({\bf D2})} and \emph{({\bf R1})}, then also \emph{({\bf R0})} is satisfied, with the parameters $ \kappa, K,\xi_0 $ depending only on $ \xi_1 $ and $ \nu $.
\end{lemma}

The part of the proof considering the exponentially decaying correlation matrices relies on the following technical result that is proven in the appendix.

\NLemma{Jensen-Poisson bound}{
Suppose $ f $ is an analytic function on the complex strip,
\bels{def of R_nu}{
\R_\nu \,:=\, \R + \cI\2(-\nu,+\nu)
\,,
}
of width $ \nu > 0$.
If $ f $ satisfies
\bels{JP: assumption on f}{
\sup_{\zeta \ins \R_\nu}\abs{f(\zeta)} \,\leq\, C_1
\qquad\text{and}\qquad
\int_0^1 \abs{f(\phi)}\2\dif \phi \,\ge\, 1
\,,
}
then for every $\eps> 0 $ there exists $ \delta > 0 $ depending only on $ \eps,\nu,C_1 $ such that 
\bels{lower bound on f}{
\absb{\setb{\phi \in [\10,1\1]:\abs{f(\phi)}\ge \delta\2}} \,\ge\, 1-\eps
\,.
}
}

\begin{Proof}[Proof of Lemma \ref{lmm:R1 and D2 imply R0}]
The non-resonance condition \eqref{non-resonance condition for A} guarantees that the $ \Lp{1}[0,1]$-norms of the row functions $ \theta \mapsto \wti{a}(\phi,\theta) $ are uniformly bounded from below. Indeed, since $ \wti{a}(\phi,\theta) \ge 0 $, we have
\bels{}{
\norm{\2\wti{a}(\phi,\genarg)}_1 
\,=\, 
\int_0^1\! \wti{a}(\phi,\theta)\2\dif \theta 
\,=\,
\sum_{j\ins\T} \nE^{\cI\12\1\pi\1j\phi}a_{j0} \,\ge\,\xi_1
\,.
}

From the exponential decay assumption ({\bf D2}) it follows that the kernel function $ \wti{a} $ has an analytic extension to the complex strip $ \R_\nu $, where $ \nu > 0 $ is the exponent from  \eqref{exponential decay of A}.
Using Lemma \ref{lmm:Jensen-Poisson bound} with $ f(\zeta) = \wti{a}(\phi,\zeta)/\xi_1 $ for a fixed $ \phi $ we see that for any $ \eps > 0 $ there exists $ \delta > 0 $ depending only on $\eps $ such that
\bels{wht-a marginal lower bound}{
\absb{ \setb{\1\theta \ins [\10\1,1] : \wti{a}(\phi,\theta) \ge \delta}} \,\ge\, 1-\eps
\,,\qquad\forall\,\phi \ins [\10\1,1\1]
\,.
}
Using \eqref{wht-a marginal lower bound} we now show that $ \wti{A} $ is a block fully indecomposable operator, i.e.,  ({\bf R0}) holds. From  \eqref{C^1-continuity of wti-a} we see that
\bels{wti-a is uniformly continuous}{
\abs{\2\wti{a}(\phi_1,\theta_1)-\2\wti{a}(\phi_2,\theta_2)} 
\,\lesssim\, 
\abs{\phi_1-\phi_2}+\abs{\theta_1-\theta_2}
\,,
}
for every $ \phi_1,\phi_2,\theta_1,\theta_2 \in [0,1] $.
Let $ K \in \N $ be so large that 
\[
\absb{\2\wti{a}(\phi_1,\theta_1)-\2\wti{a}(\phi_2,\theta_2)} \,\leq\, \frac{\delta}{2}
\,,\quad\text{provided}\quad
\abs{\phi_1-\phi_2}+\abs{\theta_1-\theta_2} \leq \frac{1}{K}
\,.
\]
Let us define a partition $ \mcl{D} = \sett{D_k}_{k=1}^K $ of $ [\10\1,1] $ and a matrix $ \brm{Z} = (Z_{ij})_{i,j=1}^K $, by
\bels{}{
D_k := \biggl[\frac{k-1}{K},\frac{k}{K}\biggr)
\qquad\text{and}\qquad
Z_{ij} := \Ind\setbb{\max_{\;(\phi,\theta)\ins D_i \times D_j\msp{-10}}\,\wti{a}(\phi,\theta) \ge \delta}
\,.
}
By the choice of $K$, we have
\bels{}{
\wti{a}(\phi,\theta) \,\ge\,\frac{\delta}{2}\,Z_{ij}\,,\qquad
(\phi,\theta) \in D_i \times D_j
\,. 
}
We will now show that $ \brm{Z} $ is fully indecomposable by proving that if there are two sets $I $ and $ J$ such that $ Z_{ij} = 0 $, for all $ i \in I $ and $ j \in J $, then
\bels{abs-I + abs-J bounded by K-1}{
\abs{\1I\1} + \abs{J} \,\leq\, K-1
\,.
}
Denoting $ D_I := \cup_{i\ins I} D_i $, we have $ \wti{a}(\phi,\theta) \leq \delta $ for $ (\phi,\theta) \in D_I \times D_J $. Thus \eqref{wht-a marginal lower bound} implies
\bels{}{
\frac{\abs{\1I\1}}{K} = \abs{D_I} \leq \eps\,,\qquad
\text{and}\qquad
\frac{\abs{J}}{K} = \abs{D_I} \leq \eps
\,.
}
Choosing $ \eps \leq 1/3 $ we see that $ \abs{I} + \abs{J} \leq (2/3)K $, and \eqref{abs-I + abs-J bounded by K-1} follows. Since $ \brm{Z} $ is a fully indecomposable matrix we see that $ \wti{A} $ is block fully indecomposable.  
\end{Proof}

\begin{lemma}[Expected decay of off-diagonal resolvent entries]
\label{lmm:Expected decay of off-diagonal resolvent entries}
If $ \brm{A} $ satisfies \emph{({\bf R0})}, in addition to \emph{({\bf D1})} or \emph{({\bf D2})}, then \eqref{decay of q_x(z)} holds with the constant $ C $ depending only on the model parameters.
\end{lemma}

\begin{Proof}
Recall from Lemma \ref{lmm:Bounded solution} that $ \wti{m}(z) $ is the bounded solution of the continuous QVE \eqref{continuous QVE}. 
We will first prove that 
\bels{def of wht-q}{
\wti{q}_x(z) \,:=\, \avg{\1e_x,\wti{m}(z)}
\,,\qquad x \in \Z
\,, 
}
satisfies
\bels{wti-q_x(z) decay estimates}{
\abs{\2\wti{q}_x(z)\1} \,\lesssim\, 
\begin{cases}
\displaystyle
\,\abs{x}^{-\1\kappa}
\quad &
\text{when \eqref{beta-summability of a_xy} holds;}
\vspace{0.1cm}
\\
\displaystyle
\nE^{-\1\nu'\abs{\1x\1}} &\text{when \eqref{exponential decay of A} holds;}
\end{cases}
\qquad x \in \Z\,,
}
where  $ e_x $ is the Fourier-basis function. 
Then we show that $ q_x(z) $ and $ \wti{q}_x(z) $ are so close that \eqref{decay of q_x(z)} holds.

Let us first assume that $ \brm{A} $ is exponentially decaying, i.e., ({\bf D2}) holds.
Let us periodically extend the kernel function $ \wti{a} : [0,1]^2 \to [0,\infty) $ from \eqref{def of wti-a} to all of $\R^2$. 
From \eqref{exponential decay of A} it follows that $ \wti{a} $ can be further analytically extended to the product of complex strips $ \R_{\nu/2}^2 $ (cf. \eqref{def of R_nu}), where $ \nu > 0 $ is the exponent from  \eqref{exponential decay of A}. 
We will now show that $ \wti{q}_x(z) $ decays exponentially in this case. To see this we consider the function $ \Gamma(z) : \R_\nu \to \C $, defined by 
\bels{def of Gamma(z;zeta)}{
\Gamma(z;\zeta) \,:=\, -\2\biggl(z+\int_0^1 \wti{a}(\zeta,\phi)\1\wti{m}(z;\phi)\2\dif \phi\biggr)^{\!-\11}
\,.
}
In particular, it follows that $ \wti{m}(z;\phi) = \Gamma(z;\phi) $ for every $ \phi \in [\10\1,1] $.
Because $ \wti{a} $ is uniformly continuous and the expression inside the parenthesis on the right hand side of \eqref{def of Gamma(z;zeta)} is bounded away from zero by a constant comparable to $ (\1\sup_z\norm{\wti{m}(z)}_\infty)^{-1} $ when $\zeta \in \R$, there exists a constant  $ \nu' < \nu $ such that $ \abs{\1\Gamma(z;\zeta)} \leq C_0 $ for $ \zeta \in \R_{\nu'} $.
Since $ \wti{a} :  \R_{\nu'}^2 \to \C $ is analytic also $ \Gamma(z) : \R_{\nu'} \to \C $ is analytic. 
For any $ x \in \Z $ we thus get by a contour deformation
\bels{q_x contour integral}{
\nE^{\12\pi\1\nu'x}\wti{q}_x(z) \,=\, \nE^{\12\pi\1\nu'x}\avg{\1e_x\1,\wti{m}(z)} 
\,&=\, \int_0^1 \!\nE^{-\1\cI\12\1\pi\1x\2(\phi\2+\2\cI\1\nu')}\,\Gamma(z;\phi)\2\dif \phi 
\\
&=\,\int_0^1 \!\nE^{-\1\cI\12\nu\1x\1\phi}\, \Gamma(z;\phi-\cI\1\nu')\2\dif \phi
\,,
}
where the integrals over the vertical line segments joining $ \pm 1 $ and $ \pm 1 - \cI\1\nu' $ cancel each other due to periodicity of the integrand in the horizontal direction. 
Since $ x \in \Z $  was arbitrary,  taking absolute values of \eqref{q_x contour integral} yields the exponential decay:
\bels{exponential decay wti-q_x(z)}{
\abs{\1\wti{q}_x(z)} 
\;\leq\;
\biggl(\sup_{\;\zeta \ins \R_{\nu'}\!}\abs{\2\Gamma(z;\zeta)}\,\biggr)\,\nE^{-2\pi\2\nu'\abs{x}} 
\;\leq\;
C_0\2\nE^{-2\pi\2\nu'\abs{x}} 
\,,
\qquad x \in \Z
\,.
}

Next we prove that ({\bf D1}) implies $ \abs{\wti{q}_x(z)} \lesssim \abs{x}^{-\kappa} $.
To this end let $ \partial $ denote the derivative w.r.t. the variable in $ [0,1] $. Using $ e_x(\phi) = \nE^{\2\cI\12\pi\1x\2\phi} $ we get for any $ k \in \N $:
\bels{x^k wti-q_x bounded by partial^km}{
\abs{x}^k \abs{\1\wti{q}_x(z)} \,=\, (2\pi)^{-k}\abs{\avg{\partial^{\1k}\msp{-2}e_x,\wti{m}(z)\1}} \,=\,(2\pi)^{-k}\abs{\avg{\1e_x,\partial^{\1k}\wti{m}(z)\1}} \leq \norm{\1\partial^{\1k}\wti{m}(z)}_\infty
\,,\quad
\forall\1x \in \Z
\,.
}
Thus, we need to show that $ \norm{\1\partial^{\1\kappa}\wti{m}(z)}_\infty \lesssim 1 $ uniformly in $ z $.
The proof is by induction on the number of derivatives of $ \wti{m} $. It is based on 
\[
\partial^{\1k}\wti{m}(z;\phi) \,=\; \partial_{\!\phi}^{\1k-1}\biggl( \wti{m}(z;\phi)^2\! \int_0^1\!\wti{m}(z;\theta)\2 \partial_\phi\1\wti{a}(\theta,\phi)\1\dif\theta\,\biggr)
\,,
\]
which follows from  \eqref{differential of wti-m},  and the following consequence of \eqref{beta-summability of a_xy}:
\bels{wht-a derivatives}{
\max_{j\1=\10}^{\kappa}\,
\sup_{\phi,\1\theta \ins [\10\1,1]} \absb{\1\partial^{\1j}_{\!\phi}\wti{a}(\phi,\theta)}
\,\lesssim\, 1
\,.
}

As the second step of the proof we show that 
\bels{q_x is bounded by wti-q_x}{
\absb{\1\wti{q}_x(z)\2-\,q_x(z)}
\;\lesssim\, 
N^{-1/2}
\,,
\qquad \abs{x}\leq N/2 +1
\,,
}
where we represent $ \T $ by integers $ x $ satisfying $\abs{x} \leq N/2+1$.
Combining \eqref{q_x is bounded by wti-q_x} with \eqref{wti-q_x(z) decay estimates} yields \eqref{decay of q_x(z)}. 
To get \eqref{q_x is bounded by wti-q_x} we use \eqref{uniform bounds of wti-m} and \eqref{m and wti-m difference} to obtain
\bea{
&\absb{\1\wti{q}_x(z)\2-\,q_x(z)}
\;\leq\;
\absbb{\,
\int_0^1 \nE^{-\cI\12\pi\1x\1\phi}\2\wti{m}(z;\phi)\1\dif\phi
- \frac{1}{N}\sum_{\theta\ins\Tc} \nE^{-\cI\12\pi\1x\1\theta}m_\theta(z)
\,}
\\
&\msp{70}\leq\;
\sum_{j=0}^{N-1} \int_{j/N}^{(j+1)/N}
\absB{\,
\nE^{-\cI\12\1\pi\1x\1\phi}\2\wti{m}(z;\phi)-\nE^{-\cI\12\1\pi\1x\1\frac{j}{N}}m_{j/N}(z)
}\,\dif \phi
\;\,\lesssim\;
\frac{1+\abs{x}}{N} \,+\, \frac{1}{\2N^{1/2}}
\,.
} 
This proves \eqref{q_x is bounded by wti-q_x} for $ \abs{x} \leq N^{1/2} $.

For $ \abs{x} \ge N^{1/2}$ we bound $ q_x = q_x(z) $ directly by using the summation of parts
\[
q_x
\,=\,
\frac{1}{N}
\sum_{j=0}^{N-1}\nE^{-\cI\12\1\pi\1x\frac{j}{N}}m_{j/N}
\,=\,
-\frac{1}{N}\sum_{j=1}^{N-2}\Bigl(m_{(j+1)/N}-m_{j/N}\Bigr)\sum_{k=0}^j\nE^{-\cI\12\1\pi\1x\frac{k}{N}}
\,+\,\Ord\Bigl(\frac{1}{N}\Bigr)
\,,
\]
where we have dropped two boundary terms of size $ \Ord(N^{-1})$. Here, $ \abs{\1m_{(j+1)/N}-m_{j/N}} \leq C/N $, while the geometric sum is $ \Ord(N/x) = \Ord(N^{1/2})$  for $ N^{1/2} \leq \abs{x} \leq N/2+1$. Thus, estimating each term in the sum over $ j $ separately shows that $ \abs{q_x(z)} \lesssim N^{-1/2} $ also in this case.
\end{Proof}

Next we show that the probability density $ \rho $ corresponding to the discrete QVE, via \eqref{def of density of states}, is also regular and supported on a single interval.

\begin{Proof}[Proof of Proposition \ref{prp:Regularity of q_x}]
Uniform boundedness of $ \brm{m} $ follows from Lemma \ref{lmm:Bounded solution}.
The other statements concerning the density $ \rho $ follow by using Theorems 1.1, 1.2, and 1.9 from \cite{AEK1}. As an input for Theorem 1.9 in \cite{AEK1}, which shows that the support of $ \rho $ is a single interval, we use 
\[
\sup_{D \1\subset\1 \Tc} \inf_{\substack{\phi_1 \ins D\\ \phi_2 \1\notin\1 D}} \sum_{\theta \in \Tc}\,\absb{\1\wht{a}_{\phi_1,\theta}-\wht{a}_{\phi_2,\theta}}
\;\leq\; 
\frac{C}{N} \,+\! \sup_{D \1\subset\1 [\10\1,1\1]} \inf_{\substack{\phi_1 \ins D\\ \phi_2 \1\notin\1 D}} 
\normb{\2\wti{a}(\phi_1,\genarg)-\2\wti{a}(\phi_2,\genarg)}_1 
\;\lesssim\; 
\frac{1}{N}
\,.
\]
The first bound follows from $ \abs{N\1\wht{a}_{\phi\1\theta}-\wti{a}(\phi,\theta)}\leq CN^{-1} $. The last bound follows from \eqref{single interval condition for wti-a}. 
The components $ q_x(z) $ inherit their analyticity trivially from $ \brm{m}(z) $ since the sum in the definition \eqref{def of q_x(z)} of $ q_x(z) $ is absolutely summable.    
\end{Proof}

\section{Proofs for local law and bulk universality}
\label{sec:Proofs for local law and bulk universality}

The following is the strongest version of the local law we prove here.

\begin{proposition}[Local law]
\label{prp:Local law}
Let $\brm{H} $ and $ \brm{A} $ be related by \eqref{EE H H = A + B}. 
Assume that $ \brm{A} $ satisfies \emph{({\bf R0})} and \emph{({\bf D1})}.
Then the conclusions of Theorem \ref{thr:Local law for Gaussian matrices with correlated entries} hold.
\end{proposition}

\begin{Proof}[Proof of Theorem \ref{thr:Local law for Gaussian matrices with correlated entries}]
If ({\bf D1}) and ({\bf R2}) are assumed, then ({\bf R0}) holds with $ \xi_0 = \xi_2 $ and $ \mcl{D} = \sett{[0,1]} $, and Proposition \ref{prp:Local law} yields the proof. 
If on the other hand, ({\bf D2}) and ({\bf R1}) are assumed, then ({\bf R0}) holds by Lemma \ref{lmm:R1 and D2 imply R0}. The proof is hence again reduced to Proposition~\ref{prp:Local law}.
\end{Proof}

\begin{Proof}[Proof of Proposition \ref{prp:Local law}]
By Lemma \ref{lmm:Fourier transform} the Fourier transform $ \wht{\brm{H}} $ of $ \brm{H} $ has the correlation structure \eqref{H:first two moments in p-space}. 
In particular, $ \wht{\brm{H}} $ is $4$-fold correlated (Definition \ref{def:4-fold correlated ensemble}). Moreover, from \eqref{H: covariance in p-space} we read off that
\[ 
\EE\,\wht{h}_{\phi\1\theta}\1\wht{h}_{-\theta,\phi} \,=\, 0
\,,\qquad\forall\,\phi,\theta \in \Tc\,,\quad \phi\neq \theta 
\,.
\]  
Hence the local law for $ 4$-fold correlated matrices, Theorem \ref{thr:Local law for 4-fold correlation}, with $ \wht{\brm{H}} $ and $ \wht{\brm{A}} $ playing the roles of $ \brm{H} $ and $\brm{S}  $, applies. In particular,  \eqref{Dep Gauss: local law for DS} follows.

In order to get \eqref{exponential off-diagonal decay of G} we use the anisotropic local law (Theorem \ref{thr:Anisotropic law}).
Indeed, fix two arbitrary elements $ x $ and $ y $ of $ \T $ and define two unit vectors $ \brm{v} $ and $ \brm{w} $ of $ \C^{\T} $ by setting
\[
\qquad
v_\phi := N^{-1/2}\nE^{\1\cI\12\pi\1x\1\phi} 
\qquad\text{and}\qquad
w_\theta := N^{-1/2}\nE^{\1\cI\12\pi\1y\1\theta} 
\,,
\qquad\forall\2\phi,\theta \in \Tc\,.
\]
From \eqref{def of F-transform} and \eqref{def of q_x(z)} we see that
\[ 
G_{xy}(z) \,=\, \brm{v} \cdot \wht{\brm{G}}(z)\1\brm{w}
\qquad\text{and}
\qquad q_{x-y}(z) \,=\, \brm{v} \cdot \mrm{diag}(\brm{m}(z))\brm{w}
\,,
\]
where $ \brm{v} \cdot \brm{w} = \sum_i \overline{\2v_i}\2w_i $. 
Thus the anisotropic local law \eqref{isoeq} implies \eqref{exponential off-diagonal decay of G}. 
The decay estimate \eqref{decay of q_x(z)} for $ q_x $ is already proven in  Lemma \ref{lmm:Expected decay of off-diagonal resolvent entries}.
\end{Proof}

Next we show that the eigenvalues of $ \brm{H} $ satisfy also the bulk universality provided the elements of $ h_{ij} $ contain a small Gaussian GOE/GUE component.

\begin{Proof}[Proof of Corollary \ref{crl:Bulk universality}]
We will show that there exists a Gaussian random matrix $ \brm{H}^{(0)} $ and a GOE/GUE matrix $ \brm{U} $ that is independent of $ \brm{H}^{(0)} $, such that the Gaussian random matrix
\bels{H decomposition}{
\brm{H} \,=\, 
\brm{H}^{(0)} + \sqrt{\frac{\1\eps\1}{2}}\2\brm{U} 
\,,
}
where $ \eps > 0 $ equals either $ \xi_2 $ or $ \xi_3 $ depending on the symmetry class, satisfies \eqref{H:first two moments in x-space}.
The matrix $ \brm{H}^{(0)} $ is such that Theorem \ref{thr:Local law for Gaussian matrices with correlated entries} is applicable since the associated correlation matrix $ \brm{A}^{\!(0)} $ satisfies ({\bf D1}) and ({\bf R2}).
In particular, the eigenvalues of $ \brm{H}^{(0)} $ satisfy the rigidity estimate (1.33) of Corollary~1.10 in \cite{AEK2}. Here we note that the corollary holds trivially for the four-fold correlated matrix as its proof depends only on the local law and not on the dependence structure of $ \brm{H}^{(0)} $.  
The bulk universality is hence proven exactly the same way as Theorem 1.15 in \cite{AEK2}, using the method of \cite{EK}.

Let us first consider the case where $ \brm{H} $ is real symmetric. 
In this case, the equations \eqref{H:first two moments in x-space} hold with $ \brm{B} = \brm{A} $. Comparing this with the GOE correlation structure,
\[
\EE\,u_{ij}u_{kl} = \frac{1}{N}
(\2\Delta_{i-k,j-l} + \Delta_{i-l,j-k}\1)
\,,
\qquad
\Delta_{xy} := \delta_{x0}\1\delta_{y0}\,,
\] 
we see that $ \brm{U} $ also satisfies \eqref{H:first two moments in x-space} with $ \brm{U} $ and $ \brm{\Delta}$ in place of $ \brm{H} $ and $ \brm{A}=\brm{B} $, respectively.
Applying Lemma \ref{lmm:Fourier transform} and using $ \wht{\Delta}_{\phi\1\theta} = N^{-1} $ we  obtain the representation, 
\bels{linear filtering in p-coords}{
\wht{h}_{\phi\1\theta} \,=\, \sqrt{N\2\wht{a}_{\phi\1\theta}}\;\wht{v}_{\phi\1\theta}
\,,
}
where $ \wht{v}_{\phi\1\theta} $ are the components of the Fourier transform of some GOE matrix $ \brm{V} $.
Using this representation we can identify the matrix $ \brm{H}^{(0)} $ in \eqref{H decomposition}. Namely, we define it in Fourier-coordinates,
\bels{}{
\wht{h}^{(0)}_{\phi\1\theta} := \sqrt{N\1\wht{a}_{\phi\1\theta}-\frac{\eps}{2}\,}\,\wht{v}^{(0)}_{\1\phi\1\theta} 
\,,
} 
where the term in the square root is bounded from below by $ \eps/2 $ by the assumption $ \wht{a}_{\phi\1\theta} \ge \eps/N $, and $ \brm{V}^{(0)} $ is a GOE random matrix that is independent of $ \brm{U} $. Since $ \brm{H}^{(0)} $ and $ \brm{U} $ are independent
\[ 
\wht{h}_{\phi\1\theta} \,:=\, \wht{h}^{(0)}_{\phi\1\theta} \2+ \sqrt{\frac{\1\eps\1}{2}}\,\wht{u}_{\phi\1\theta}
\,,
\] 
satisfy \eqref{H:first two moments in p-space}. This immediately yields \eqref{H:first two moments in x-space} for the matrix \eqref{H decomposition}.

Next we consider the case where $ \brm{H} $ is complex self-adjoint.
First we remark that for a given pair of hermitian matrices $ (\brm{A},\brm{B}) $ there exists a random matrix $ \brm{H} $ satisfying \eqref{H:first two moments in x-space} if and only if the following hold in the Fourier-space: 
\bels{A vs B in p-coords}{
\wht{a}_{\1\phi\1\theta} \,\ge\, 0\,,\qquad
\wht{b}_{-\phi,-\theta} \,=\, \wht{b}_{\1\phi\1\theta}
\,,
\qquad
\text{and}
\qquad  
\absb{\2\wht{b}_{\phi\1\theta}} \,\leq\, \sqrt{\2\wht{a}_{\phi\1\theta}\,\wht{a}_{-\phi,-\theta}\!}
\;, 
\qquad
\forall\2\phi,\theta \in \Tc
\,.
}
The necessity of these conditions follow from $ \wht{a}_{\phi\1\theta} = \EE\,\abs{\wht{h}_{\phi\1\theta}}^2$, $ \wht{b}_{\phi\1\theta} = \EE\,\wht{h}_{\phi\2\theta}\wht{h}_{-\phi,-\theta} $ (cf. \eqref{H: covariance in p-space}), and the Cauchy-Schwartz inequality,
\[ 
\wht{b}_{\phi\1\theta} 
\,\leq\, 
\sqrt{\EE\,\abs{\wht{h}_{\phi\1\theta}}^2}\,\sqrt{\EE\,\abs{\wht{h}_{-\phi,-\theta}}^2\,}
\,\leq\,
\sqrt{\,\wht{a}_{\1\phi\1\theta}\,\wht{a}_{-\phi,-\theta}\!} 
\;,\qquad
\phi \neq -\theta
\,.
\] 
If $ \theta = - \phi $, then the identity holds by definition.
In order to see that \eqref{A vs B in p-coords} is also a sufficient condition for there to exists a random matrix $ \brm{H} $ satisfying \eqref{H:first two moments in x-space} we consider a fixed index pair $ (\phi,\theta) \in \Tc^2 $. From the hermitian symmetry and Lemma \ref{lmm:Fourier transform} it follows that the two elements $ \wht{h}_{\phi,\theta} $ and $ \wht{h}_{-\phi,-\theta} $ determine the four entries of $ \wht{\brm{H}} $ that may depend on $ \wht{h}_{\phi,\theta} $. 
It is now easily checked that a given $ 4 \times 4$-real matrix $ \brm{\Gamma} $ can be a correlation matrix of the real random vector 
\[ 
\brm{x} \,:=\, \bigl(\2\Re\,\wht{h}_{\phi,\theta}\1,\2\Re\,\wht{h}_{-\phi,-\theta}\1,\2\Im\,\wht{h}_{\phi,\theta}\1,\2\Im\,\wht{h}_{-\phi,-\theta}\bigr)
\,,
\] 
if and only if it is positive-semidefinite. A simple computation reveals that $ \brm{\Gamma} $ is positive semi-definite if and only if the third condition of \eqref{A vs B in p-coords} holds. 

Assume now that $ \brm{A} $ and $ \brm{B} $ satisfy \eqref{Hermitian case: A vs B in p-coords} for some $ \eps > 0 $. Let us define $ \brm{A}^{\!(0)} $, by  
\[
\wht{a}^{(0)}_{\phi\1\theta} \,:=\, \wht{a}_{\phi\1\theta} \2-\2 \frac{\eps}{2}
\,,\qquad\forall\,\phi,\theta \in \Tc
\,.
\] 
From \eqref{Hermitian case: A vs B in p-coords} we see that $ \wht{a}_{\phi\1\theta}^{(0)} \ge \eps/2 $. Since $ (\wht{\brm{A}}^{\!(0)},\wht{\brm{B}}) $ satisfies \eqref{A vs B in p-coords}, with $ \wht{\brm{A}}^{(0)} $ in place of $ \wht{\brm{A}} $, there exists a random matrix $ \wht{\brm{H}}^{(0)} $ such that \eqref{H:first two moments in p-space} holds with $ \wht{\brm{H}}^{(0)} $ and $ \wht{\brm{A}}^{(0)} $ in place of $ \wht{\brm{H}} $ and $ \wht{\brm{A}} $, respectively.
Let $ \brm{U} $ be a GUE matrix so that \eqref{H:first two moments in x-space} holds with $ (\brm{U},\brm{\Delta},\brm{0}) $ in place of  $ (\brm{H},\brm{A},\brm{B}) $. Since  $ \wht{\Delta}_{\phi\1\theta} = N^{-1} $ we see that the Fourier transform of the random matrix \eqref{H decomposition} satisfies \eqref{H:first two moments in p-space} provided we choose $ \brm{U} $ to be independent of $ \brm{H}^{(0)} $. This is equivalent to \eqref{H:first two moments in x-space} and the proof is complete.   
\end{Proof}

From the proof of Corollary \ref{crl:Bulk universality} we read off the convolution representation for symmetric translation invariant random matrices.
\begin{Proof}[Proof of Lemma \ref{lmm:Linear filtering}]
The assumption \eqref{Bochner in x-space} implies $ \wht{a}_{\phi\1\theta} \ge 0 $.
This guarantees that $ \wht{H} $ defined through \eqref{linear filtering in p-coords} is self-adjoint.
Expressing \eqref{linear filtering in p-coords} in the original coordinates yields the representation \eqref{h_ij as convolution}.
\end{Proof}

\begin{Proof}[Proof of Corollary \ref{crl:Delocalization of eigenvectors}]
Given the anisotropic local law \eqref{isoeq} and the uniform boundedness (Lemma \ref{lmm:Bounded solution}) of the solution $ \brm{m} $ of the QVE \eqref{QVE in Tc}, the delocalization of the eigenvalues is proven exactly the same way as Corollary 1.13 in \cite{AEK2}.
\end{Proof}

\appendix 
\section{Appendix}

\begin{Proof}[Proof of Lemma \ref{lmm:Jensen-Poisson bound}]
Let $ \KK $ be an open and simply connected set with a smooth boundary, such that 
\bels{criteria for KK boundary}{
[\10\1,1\1] \subset \partial\1 \KK\,,
\qquad\text{and}\qquad
(\10\1,1\1) + \cI\1(\10,\tsfrac{2}{3}\nu)
\subset \KK \subset (-1\1,2\1) + \cI\2(\10,\tsfrac{2}{3}\nu)
\,.
}
Since $ \KK $ is in $ \R_{2\1\nu/3} $ and $ f $ is bounded and analytic on $ \R_\nu $, the assumption \eqref{JP: assumption on f} implies  %
\bels{continuity on R_nu}{
\abs{f(\xi)-f(\zeta)} &\leq C_2\,\abs{\xi-\zeta}\,,
\qquad\forall\1\xi,\zeta \in \overline{\KK}
\,,
} 
where $ C_2 < \infty $ does not depend on $ f $.
From the first inclusion of \eqref{criteria for KK boundary} it follows that 
\bels{JP - step 0}{
\absb{\setb{\phi \in [\10,1\1]:\abs{f(\phi)} < \delta}} 
\;\leq\;
\absb{\setb{\zeta \in \partial\KK : \abs{f(\zeta)}< \delta}} 
\,,
}
where $ \abs{A} $ denotes the Lebesgue measure of $ A \subset \R $.
We will prove \eqref{lower bound on f} by estimating the size of the set on the right.

Let us denote the complex unit disk by $ \DD := \sett{\zeta \in \C : \abs{\zeta} < 1} $, and let $ \zeta_0 \in \KK  $ be arbitrary. By the Riemann mapping theorem there exists a bi-holomorphic conformal map $ \Phi_{\zeta_0} : \C \to \C $ satisfying
\bels{defining properties of Phi_zeta0}{
\Phi_{\zeta_0}(\DD) = \KK
\qquad
\text{and}
\qquad
\Phi_{\zeta_0}(0) = \zeta_0\,. 
}
Since the simple connected sets $ \DD $ and $ \KK $ have smooth boundaries the conformal map $ \Phi_{\zeta_0} $ extends to the boundary, such that $ \Phi(\partial\DD) = \partial \KK $, with uniformly bounded derivatives. In particular, we have
\bels{bounds on derivative of Phi}{
\frac{1}{C_3(\zeta_0)} 
\,\leq\, 
\abs{\1\partial\1 \Phi_{\msp{-2}\zeta_0}(\zeta)}
\,\leq\, 
C_3(\zeta_0)\,,
\qquad\zeta \in \DD\,,
}
with the constant $ C_3(\zeta_0) < \infty $ independent of $ f $, in fact it depends only on $ \zeta_0 $ through the distance $ \dist(\1\zeta_0,\partial \KK\1) $.
From the second estimate of \eqref{JP: assumption on f} we know that there are points on $ [0,1] \subset \partial \KK $ where $ \abs{f} \ge 1 $. Hence using the continuity  \eqref{continuity on R_nu} we may choose $ \zeta_0 \in \KK $ such that 
\bels{criteria for zeta_0}{ 
\abs{f(\zeta_0)} \,\ge\, \frac{1}{2}
\qquad\text{and}\qquad
\dist(\1\zeta_0\1,\partial\KK\1) 
\,\ge\, 
\min\setbb{\frac{1}{2\1C_2\!}\2,\,\frac{\nu\!}{3}}
\,.
} 

Let $ \log_+ $ and $ \log_- $ be the positive and negative parts of the logarithm, respectively, so that $ \log \tau = \log_+\tau -\log_-\tau $, for $ \tau > 0$. 
Using Chebyshev's inequality we get
\[
\absb{\setb{\zeta \in \partial\KK : \abs{f(\zeta)}< \delta}} 
\,\leq\,
\frac{1}{\log_- \delta\,}\!
\int_{\partial\KK} \log_-\abs{f(\zeta)}\,\abs{\dif\zeta}
\,.
\]
By parametrizing the boundary of $ \KK $ using the conformal map $ \Phi_{\zeta_0} $ we get
\[
\absb{\setb{\zeta \in \partial\KK : \abs{f(\zeta)}< \delta}} 
\,\leq\,
\frac{1}{\log_- \delta\,}\!
\int_0^{2\pi} \!\log_-\abs{f(\Phi_{\zeta_0}(\nE^{\cI\1\tau}))}\,\abs{\1\partial\Phi_{\zeta_0}(\nE^{\cI\1\tau})}\2\dif \tau
\,.
\]
Using \eqref{bounds on derivative of Phi} to bound the derivative and writing $ \wti{f} := f \circ \Phi_{\zeta_0} $ we get 
\bels{JP: step 2}{
\absb{\setb{\zeta \in \partial\KK : \abs{f(\zeta)}< \delta}} 
\,\leq\,
\frac{C_3(\zeta_0)}{\log_- \delta}
\int_0^{2\pi}\! \log_-\abs{\2\wti{f}(\nE^{\cI\1\tau})}\2\dif \tau
\,.
}
We will now bound the last integral using the Jensen-Poisson formula,
\[
\log\, \abs{\wti{f}(0)} \,=\, \frac{1}{2\pi}\int_0^{2\pi} \!\log\, \abs{\1\wti{f}(\nE^{\1\cI\1\tau})}\2\dif\tau \,- \sum_{j=1}^n \log \frac{1}{\abs{\alpha_j}}
\,,
\]
where $ \alpha_j $'s are the zeros of $ \wti{f} $ in the unit disk $ \DD $. The last sum is always non-negative since $ \abs{\alpha_i} \leq 1 $ and can be dropped. By splitting the integral into positive and negative parts we obtain an estimate for the integral on the right hand side of \eqref{JP: step 2}
\bea{
\int_0^{2\pi} \!\log_- \abs{\wti{f}(\nE^{\1\cI\1\tau})}\2\dif\tau
\,
&\leq\, 
2\1\pi \log \frac{1}{\abs{\wti{f}(0)}} 
 + \int_0^{2\pi} \!\log_+ \abs{\wti{f}(\nE^{\1\cI\1\tau})}\2\dif\tau
\\
&\leq\, 
2\1\pi \log 2\,+ 2\1\pi\2 \log \sup_{\omega \ins \DD} \abs{\wti{f}(\omega)} 
\\
&\leq\,
2\1\pi \log 2\1C_1
\,,
} 
where we have used \eqref{criteria for zeta_0} to get the second inequality. For the last bound we have used $ \abs{\1\wti{f}(\omega)} = \abs{f(\Phi_{\zeta_0}(\omega))} \leq C_1 $.
Plugging this into \eqref{JP: step 2} and recalling \eqref{JP - step 0} we get
\[
\absb{\setb{\phi \in [\10,1\1]:\abs{f(\phi)} < \delta}} 
\,\leq\,
\frac{2\1\pi\1C_3(\zeta_0)\log 2C_1}{\log (1/\delta)}
\,.
\]
This finishes the proof as $ C_3(\zeta_0) $ and $ C_1 $ are independent of $ \delta $. 
\end{Proof}


\end{document}